\documentclass[twoside,12pt,a4wide]{article}
\usepackage{amsfonts}
\usepackage{amssymb}
\usepackage{a4wide}
\usepackage{amsthm}
\usepackage{amsmath}
\usepackage[all]{xy}
\usepackage{enumerate}
\usepackage{fancyhdr}

\newtheorem{theorem}{Theorem}[section]
\newtheorem{proposition}[theorem]{Proposition}
\newtheorem{corollary}[theorem]{Corollary}
\newtheorem{lemma}[theorem]{Lemma}
\theoremstyle{definition}
\newtheorem*{notation}{Notation}
\newtheorem*{Beweis}{Proof}
\newtheorem{definition}[theorem]{Definition}
\newtheorem{punto}[theorem]{}
\theoremstyle{remark}
\newtheorem{remark}[theorem]{Remark}
\newtheorem{ex}[theorem]{Example}
\newtheorem{exs}[theorem]{Examples}
\newtheorem{remarks}[theorem]{Remarks}
\CompileMatrices
\swapnumbers
\CompileMatrices
\input{tcilatex}
\begin{document}

\title{Semicorings and Semicomodules}
\author{\textbf{Jawad Y. Abuhlail}\thanks{%
The author would like to acknowledge the support provided by the Deanship of
Scientific Research (DSR) at King Fahd University of Petroleum $\&$ Minerals
(KFUPM) for funding this work through project No. IN080400.} \\
Department of Mathematics and Statistics\\
Box $\#\;$5046, KFUPM, 31261 Dhahran, KSA\\
abuhlail@kfupm.edu.sa}
\date{\today }
\maketitle

\begin{abstract}
In this paper, we introduce and investigate \emph{semicorings} over
associative semirings and their categories of \emph{semicomodules.} Our
results generalize old and recent results on corings over rings and their
categories of comodules. The generalization is \emph{not} straightforward
and even subtle at some places due to the nature of the base category of
commutative monoids which is neither Abelian (not even additive) nor
homological, and has no non-zero injective objects. To overcome these and
other difficulties, a combination of methods and techniques from
categorical, homological and universal algebra is used including a new
notion of exact sequences of semimodules over semirings.
\end{abstract}

\section*{Introduction}

\qquad Coalgebraic structures in general, and categories of comodules for
\emph{comonads} in particular, are gaining recently increasing interest \cite%
{Wis2012}. Although comonads can be defined in arbitrary categories, nice
properties are usually obtained in case the comonad under consideration is
isomorphic to $-\bullet C$ (or $C\bullet -$) for some \emph{comonoid} $C$
\cite{Por2006} in a monoidal category $(\mathbb{V},\bullet ,\mathbf{I})$
\emph{acting} on some category $\mathbb{X}$ in a nice way \cite{MW}.

However, it can be noticed that the most extensively studied \emph{concrete
examples} are -- so far -- the categories of \emph{coacts }(usually called
\emph{coalgebras}) of an endo-functor $F$ on the category $\mathbf{Set}$ of
sets (motivated by applications in theoretical computer science \cite%
{Gum1999} and universal (co)algebra \cite{AP2003}) and categories of
comodules for a coring over an associate algebra \cite{BW2003}, \cite%
{Brz2009}, \emph{i.e.} a comonoid in the monoidal category $(_{A}\mathbf{M}%
_{A},\otimes _{A},A)$ of $(A,A)$-bimodules, for some associative algebra $A,$
which acts on the category $\mathbf{M}_{A}$ ($_{A}\mathbf{M}$) of right
(left) $A$-modules in the obvious way \cite[p. 229]{Por2008}.

The main goal of this paper is to investigate categories of \emph{%
semicomodules} for a \emph{semicoring}, which can be seen as comodules of
comonads associated to a comonoid in the monoidal category $(_{A}\mathbb{S}%
_{A},\otimes _{A},A)$ of $(A,A)$-bisemimodules over an associative
semialgebra $A$ with the natural tensor product $-\otimes _{A}-$ \cite%
{Kat1997}. This does not only add a \emph{new concrete example} where the
general theory of comonads and comonoids applies, but also provides an
interesting context where a combination of techniques and methods from
categorical algebra, homological algebra and universal algebra applies
naturally and harmonically.

Semicorings over semirings are of particular importance for theoretical and
practical reasons: On one hand, and in contrast to categories of modules
over a ring, categories of semimodules over a semiring are \emph{not} so
nice, as (in general):

\begin{itemize}
\item they are not \emph{Grothendieck}: a category of semimodules over a
semiring does not necessarily have \emph{enough injectives} ($0$ is the only
injective object in the category $\mathbf{AbMonoid}\simeq \mathbb{S}_{%
\mathbb{N}_{0}},$ the category of semimodules over the semiring $\mathbb{N}%
_{0}$ of non-negative integers \cite[17.21]{Gol1999}).

\item they are not \emph{Abelian}: a semimodule does not necessarily posses
a projective presentation (objects with a projective presentation are called
\emph{normal} \cite{Tak1983}); moreover, one cannot make free use of \emph{%
some} nice properties of adjoint functors between Abelian categories (\emph{%
e.g.} \cite[Proposition 6.28]{Fai1973}); see Remark \ref{pres-inj}.

\item they are not \emph{additive}; the hom sets are monoids which are not
necessarily groups.

\item they are not \emph{Puppe-exact} \cite{Pup1962}; for example, the
notion of \emph{exact sequences of semimodules} is subtle; this led to
several different notions of exactness for sequences of semimodules over a
semiring \cite{Abu-b} (all of which coincide for modules over rings).

\item they are not \emph{homological} since they are not \emph{protomodular}
\cite{BB2004}: several basic homological lemmas do not apply (\emph{e.g.}
the short five lemma). Moreover, epimorphisms are not necessarily surjective
and subsemimodules are not necessarily kernels.

\item several flatness, projectivity and injectivity properties which are
equivalent for modules over rings are apparently different for semimodules
over semirings \cite{Abu-a}.

\item Some notions cannot be easily checked; for example, to prove that a
given $S$-semimodule is \emph{flat}, one has to show that $M\otimes _{S}-:$ $%
_{S}\mathbb{S}\longrightarrow $ $\mathbf{AbMonoid}$ preserves all pullbacks
\cite{Kat2004} (or all equalizers) and not only the monomorphisms.
\end{itemize}

This has the impact that generalizing results on corings (comodules) to
semicorings (semicomodules) is neither trivial nor straightforward as the
first impression might be. We could overcome some of the difficulties
mentioned above by introducing a new notion for exact sequence of
semimodules over semirings which we used to prove \emph{restricted versions}
of the short five lemma \cite{Abu-a} (the nine lemma and the snake lemma
\cite{Abu-b}). We also introduced and used suitable notions of flatness,
projectivity and injectivity for semimodules. Moreover, we made use of
recent developments in the theory of comonads \cite{BBW2009} especially
those associated to comonoids in monoidal categories \cite{Por2006}, \cite%
{Por2008b}.

On the other hand, semirings and semimodules proved to have a wide spectrum
of significant applications in several aspects of mathematics like
optimization theory \cite{C-G1979}, tropical geometry \cite{R-GST2005},
idempotent analysis \cite{KM1997}, physics \cite{Gun1998}, theoretical
computer science (\emph{e.g.} Automata Theory \cite{Eil1974}, \cite{Eil1976}%
) and many more \cite{Gol1999}. Moreover, corings over rings showed to have
important applications in areas like noncommutative ring theory, category
theory, Hopf algebras, differential graded algebras, and noncommutative
geometry \cite{Brz2009}. This suggests that investigating semicorings and
semicomodules will open the door for many new applications in the future
(see \cite{Wor2012} for recent applications to Automata Theory).

Before proceeding, we mention that in many (relatively old) papers,
researchers used the so called Takahashi's \emph{tensor-like} product, which
we denote by $-\boxtimes _{A}-$ \cite{Tak1982} (see also \cite{Gol1999}).
This product has the defect that, for a semialgebra $A,$ the category $(_{A}%
\mathbb{S}_{A},\boxtimes _{A},A)$ of $(A,A)$-bisemimodules is \emph{not}
monoidal. The author \cite{Abu2013} introduced a notion of \emph{semiunital
semimonoidal categories }with prototype $(_{A}\mathbb{S}_{A},\boxtimes
_{A},A)$ with $A$ as a \emph{semiunit}; he also presented a notion of \emph{%
semicounital semicomonoids} in such categories with \emph{semicounital }$A$-%
\emph{semicorings} as an illustrating example. The relation between the two
products has been clarified in \cite{Abu-a}.

This paper is organized as follows: after this introduction, and for
convenience of the readers not familiar with semirings and semimodules, we
include in Section 1 some basic definitions, properties and some results
(mostly without proof) related to such algebraic structures. In Section 2,
we introduce the notion of a semicoring $\mathcal{C}$ over a semialgebra $A$
and study basic properties of the category $\mathbb{S}^{\mathcal{C}}$ of
right $\mathcal{C}$-semicomodules. We present a reconstruction result in
Theorem \ref{comonad-C}. Moreover, we apply results of Porst et al. (\emph{%
e.g.} \cite{AP2003}, \cite{Por2006}, \cite{Por2008}, \cite{Por2008b}) to
obtain a generalization of the Fundamental Theorem of Coalgebras over fields
to semicoalgebras over commutative semirings in Proposition \ref{semicoal}
and to semicorings over arbitrary rings in Proposition \ref{coring}. Let $%
_{A}\mathcal{C}$ be a semicoring over the semialgebra $A.$ We apply results
of Porst et al. to $\mathbb{S}^{\mathcal{C}}$ (Theorem \ref{Porst}). In
Section 3, we introduce and investigate the category $\mathrm{Rat}^{\mathcal{%
C}}(\mathbb{S}_{\mathcal{A}})$ of $\mathcal{C}$-rational right $\mathcal{A}$%
-semimodules, where $\mathcal{A}$ is an $A$-semiring with a morphism of $A$%
-semirings $\kappa :\mathcal{A}\longrightarrow $ $^{\ast }\mathcal{C}$ ($:=%
\mathrm{Hom}_{A}(\mathcal{C},A),$ the left dual ring of $\mathcal{C}$) and $%
P=(\mathcal{A},\mathcal{C})$ is a \emph{left }$\alpha $\emph{-pairing}. In
this case, we prove that $\mathrm{Rat}^{\mathcal{C}}(\mathbb{S}_{\mathcal{A}%
})\simeq \mathbb{S}^{\mathcal{C}}$ (Theorem \ref{equal}) extending our main
result in \cite{Abu2003} on the category of right comodules for a coring
over an associative algebra (see also \cite{BW2003}). Moreover, and assuming
a uniformity condition on $_{A}\mathcal{C},$ we show in Theorem \ref%
{eq-alpha} that $\mathbb{S}^{\mathcal{C}}=\sigma \lbrack \mathcal{C}_{^{\ast
}\mathcal{C}}]$ (the \emph{Wisbauer category} of $\mathcal{C}$-subgenerated
right $^{\ast }\mathcal{C}$-semimodules) if and only if $_{A}\mathcal{C}$ is
a \emph{mono-flat }$\alpha $\emph{-semimodule} and $\mathbb{S}^{\mathcal{C}}$
is closed under $^{\ast }\mathcal{C}$-subsemimodules. Under some suitable
conditions on the semialgebra $A$ and a uniformity condition $_{A}\mathcal{C}%
,$ we prove in Theorem \ref{summary} that $\mathbb{S}^{\mathcal{C}}=\mathbb{S%
}_{^{\ast }\mathcal{C}}$ if and only if $_{A}\mathcal{C}$ is finitely
generated projective and $\mathbb{S}^{\mathcal{C}}$ is closed under $^{\ast }%
\mathcal{C}$-subsemimodules.

\section{Preliminaries}

\begin{punto}
Let $(G,+)$ be an Abelian semigroup. We say that $G$ is \emph{cancellative}
iff%
\begin{equation*}
g+g^{\prime }=g+g^{\prime \prime }\Rightarrow g^{\prime }=g^{\prime \prime }
\end{equation*}%
for all $g,g^{\prime },g^{\prime \prime }\in G.$ For say that a subsemigroup
$L\leq G$ is \emph{subtractive} iff $L=\overline{L},$ where%
\begin{equation*}
\overline{L}=\{g\in G\mid g+l^{\prime }=l^{\prime \prime }\text{ for some }%
l^{\prime },l^{\prime \prime }\in L\}.
\end{equation*}%
We say that $G$ is \emph{completely subtractive} iff every subsemigroup $%
L\leq G$ is subtractive. A morphism of semigroups $f:G\longrightarrow
G^{\prime }$ is said to be \emph{subtractive} iff $f(G)\leq G^{\prime }$ is
a subtractive subsemigroup.
\end{punto}

\begin{remark}
\label{f-1}Let $f:G\longrightarrow G^{\prime }$ be a morphism of Abelian
monoids. If $L^{\prime }\leq G^{\prime }$ is a subtractive submonoid, then $%
f^{-1}(L)\leq G$ is a subtractive submonoid.
\end{remark}

\subsection*{Semirings and Semimodules}

\qquad In this section, we present some basic definitions and results on
semirings and semimodules. Our main reference is \cite{Gol1999}; however, we
use a different notion of \emph{exact sequences} of semimodules introduced
in \cite{Abu-b}. With $\mathbf{AbMonoid},$ we denote the category of Abelian
monoids.

\begin{punto}
A \emph{semiring} is a monoid in the monoidal category $(\mathbf{AbMonoid}%
,\otimes ,\mathbb{N}_{0})$ of Abelian monoids, or roughly speaking a ring
not necessarily with subtraction, \emph{i.e.} a non-empty set $S$ with two
binary operations \textquotedblleft $+$\textquotedblright\ and
\textquotedblleft $\cdot $\textquotedblright\ such that $(S,+,0)$ is an
Abelian monoid and $(S,\cdot ,1)$ is a monoid such that%
\begin{equation*}
s\cdot (s_{1}+s_{2})=s\cdot s_{1}+s\cdot s_{2},\text{ }(s_{1}+s_{2})\cdot
s=s_{1}\cdot s+s_{2}\cdot s\text{ and }s\cdot 0=0=0\cdot s\text{ for all }%
s,s_{1},s_{2}\in S.
\end{equation*}%
A \emph{morphism of semirings} $f:S\longrightarrow T$ is a map such that $%
f:(S,+_{S},0_{S})\longrightarrow (T,+_{T},0_{T})$ and $f:(S,\cdot
_{S},1_{S})\longrightarrow (T,\cdot _{T},1_{T})$ are morphisms of monoids.
We say that the semiring $S$ is \emph{commutative }(\emph{cancellative}) iff
$(S,\cdot )$ is commutative ($(S,+)$ is cancellative). If $S$ is a
commutative ring and $\eta :S\longrightarrow A$ is a morphism of semirings,
then we call $A$ an (associative) $S$\emph{-semialgebra}.
\end{punto}

\begin{punto}
Let $S$ be a semiring. A \emph{right }$S$\emph{-semimodule} $M$ is roughly
speaking a right $S$-module not necessarily with subtraction (\emph{i.e.}$%
(M,+_{M},0_{M})$ is an Abelian monoid rather than a group) for which%
\begin{equation*}
m0_{S}=0_{M}=0_{M}s\text{ for all }m\in M\text{ and }s\in S.
\end{equation*}%
The category of right (left) $S$-semimodules and $S$\emph{-linear maps},
which respect addition and scalar multiplication, is denoted by $\mathbb{S}%
_{S}$ ($_{S}\mathbb{S}$). A right (left) $S$-semimodule $M$ is said to be
\emph{cancellative} (\emph{completely subtractive}) iff $(M,+)$ is
cancellative (every $S$-subsemimodule $L\leq _{S}M$ is subtractive). With $%
\mathbb{CS}_{S}\hookrightarrow \mathbb{S}_{S}$ ($_{S}\mathbb{CS}%
\hookrightarrow _{S}\mathbb{S}$), we denote the \emph{full} subcategory of
cancellative right (left) $S$-semimodules. For semirings $S$ and $T,$ an
object in the category $_{S}\mathbb{S}_{T}$ of $(S,T)$\emph{-bisemimodules}
is a left $S$-semimodules $_{S}M$ which is also right $T$-semimodule $M_{T}$
with $s(mt)=(sm)t$ for all $s\in S,$ $m\in M$ and $t\in T;$ the arrows are
the $S$-linear $T$-linear maps, and its \emph{full} subcategory of
cancellative $(S,T)$-bisemimodules is denoted by $_{S}\mathbb{CS}_{T}.$ We
call $_{S}M_{T}$ a \emph{completely subtractive }$(S,T)$\emph{-bisemimodule}
iff every $(S,T)$-subbisemimodule $L\leq _{(S,T)}M$ is subtractive.
\end{punto}

\begin{exs}
\begin{enumerate}
\item Every ring is a semiring.

\item The set $\mathbb{N}_{0}$ of non-negative integers is a \emph{%
cancellative} semiring with the usual addition and multiplication. The
category of (cancellative) $\mathbb{N}_{0}$-semimodules is isomorphic to the
category of (cancellative) Abelian monoids.

\item An example of a semiring due to Dedekind \cite{Ded1894} is $(\mathrm{%
Ideal}(R),\cup ,\cap ),$ where $R$ is a ring and $\mathrm{Ideal}(R)$ is the
set of ideals of $R.$ More generally, any distributive complete lattice is a
semiring.

\item $\mathbb{R}_{\max }:=(\mathbb{R}\cup \{-\infty \},\max ,+)$ and $%
\mathbb{R}_{\min }:=(\mathbb{R}\cup \{\infty \},\min ,+)$ are \emph{%
semifields} (every non-zero element has a multiplicative inverse) \cite%
{C-G1979}.

\item $\mathbf{B}=(\{0,1\},+,\cdot )$ is a \emph{semi-field,} where $%
1+1=1\neq 0$ \cite[p. 7]{Gol1999}, called the \emph{Boolean semiring}.
\end{enumerate}
\end{exs}

\begin{punto}
Let $M$ be a right $S$-semimodule. Each $L\leq _{S}M$ defines two $S$\emph{%
-congruence relations} \cite{Gol1999} on $M,$ namely $\equiv _{L}$ and $%
[\equiv ]_{L},$ where%
\begin{equation*}
\begin{tabular}{l}
$m_{1}\text{ }\equiv _{L}\text{ }m_{2}\text{ iff }m_{1}+l_{1}=m_{2}+l_{2}%
\text{ for some }l_{1},l_{2}\in L;$ \\
$m_{1}\text{ }\left[ \equiv \right] _{L}\text{ }m_{2}\text{ iff }%
m_{1}+l_{1}+m^{\prime }=m_{2}+l_{2}+m^{\prime }\text{ for some }%
l_{1},l_{2}\in L\text{ and }m^{\prime }\in M.$%
\end{tabular}%
\end{equation*}%
We define the \emph{quotient semimodule} $M/L:=M/\equiv _{L};$ notice that $%
M/[\equiv ]_{L}$ is cancellative for each $L\leq _{S}M.$ In fact, we have a
functor%
\begin{equation*}
\mathfrak{c}(-):\mathbb{S}_{S}\longrightarrow \mathbb{CS}_{S},M%
\longrightarrow M/[\equiv ]_{\{0\}}.
\end{equation*}
\end{punto}

\begin{remark}
\label{tensor}The \emph{tensor product} $-\otimes _{S}-$ of semimodules we
adopt is that in the sense of \cite{Kat1997} (see also \cite{Kat2004}, \cite%
{KN2011}) and \textbf{not} in the \emph{tensor-like product} introduced by
Takahashi \cite{Tak1982} which we denote by $-\boxtimes _{S}-.$ As clarified
in \cite{Abu-a}, for every right (left) $S$-semimodule $M,$ we have an
isomorphism of Abelian monoids $M\otimes _{S}S\overset{\vartheta _{M}^{r}}{%
\simeq }M$ ($S\otimes _{S}M\overset{\vartheta _{M}^{l}}{\simeq }M$) while%
\begin{equation*}
M\boxtimes _{S}S\simeq \mathfrak{c}(M\otimes _{S}S)\simeq \mathfrak{c}(M)%
\text{ (}S\boxtimes _{S}M\simeq \mathfrak{c}(S\otimes _{S}M)\simeq \mathfrak{%
c}(M)\text{).}
\end{equation*}
\end{remark}

\begin{definition}
Let $\mathfrak{A}$ and $\mathfrak{B}$ be categories and $F:\mathfrak{A}%
\longrightarrow \mathfrak{B}$ a covariant functor.

\begin{enumerate}
\item $F$ \emph{preserves limits} iff for every diagram $D:\mathfrak{I}%
\longrightarrow \mathfrak{A}$ we have:
\begin{equation*}
(L\overset{\ell _{i}}{\longrightarrow }D_{i})_{i\in \mathrm{Obj}(\mathfrak{I}%
)}\text{ is a limit of }D\Rightarrow (F(L)\overset{\ell _{i}}{%
\longrightarrow }F(D_{i}))_{i\in \mathrm{Obj}(\mathfrak{I})}\text{ is a
limit of }F\circ D;
\end{equation*}

\item $F$ \emph{create limits} iff for every diagram $D:\mathfrak{I}%
\longrightarrow \mathfrak{A}$ and every limit $\mathcal{L}^{\prime
}=(L^{\prime }\overset{f_{i}^{\prime }}{\longrightarrow }D_{i}^{\prime
})_{i\in \mathrm{Obj}(\mathfrak{I})}$ of $F\circ D$ in $\mathfrak{B},$ if
any, there exists a unique cone $\mathcal{S}=(L\overset{f_{i}}{%
\longrightarrow }D_{i})_{i\in \mathrm{Obj}(\mathfrak{I})}$ in $\mathfrak{A}$
with%
\begin{equation*}
(F(L)\overset{F(f_{i})}{\longrightarrow }F(D_{i}))_{i\in \mathrm{Obj}(%
\mathfrak{I})}=(L^{\prime }\overset{f_{i}^{\prime }}{\longrightarrow }%
D_{i}^{\prime })_{i\in \mathrm{Obj}(\mathfrak{I})}\text{ and }\mathcal{S}%
\text{ is a limit of }D\text{ in }\mathfrak{A}.
\end{equation*}%
Dually, one defines the functor which preserving (creating) colimits.
\end{enumerate}
\end{definition}

\begin{lemma}
\label{adj-lim}\emph{(cf. \cite[Proposition 3.2.2]{Bor1994a})} Let $%
\mathfrak{A},\mathfrak{B}$ be arbitrary categories and $\mathfrak{B}\overset{%
F}{\longrightarrow }\mathfrak{B}\overset{G}{\longrightarrow }\mathfrak{C}$
be functors such that $(F,G)$ is an adjoint pair.

\begin{enumerate}
\item $F$ preserves all colimits which turn out to exist in $\mathfrak{A}.$

\item $G$ preserves all limits which turn out to exist in $\mathfrak{B}.$
\end{enumerate}
\end{lemma}

\begin{definition}
An object $G$ in a cocomplete category $\mathfrak{A}$ is said to be a (\emph{%
regular})\emph{\ generator }iff for every $X\in \mathfrak{A},$ there exists
a canonical (\emph{regular})\emph{\ epimorphism} $f_{X}:\dbigsqcup\limits_{f%
\in \mathfrak{A}(G,X)}G\longrightarrow X$ \cite[p. 199]{BW2005} (see also
\cite{Kel2005}, \cite{Ver}); recall that an arrow in $\mathfrak{A}$ is said
to be a \emph{regular epimorphism} iff it is a coequalizer (of its \emph{%
kernel pair}).
\end{definition}

The following result collects some properties of the categories of
semimodules over a semiring (cf. \cite{Abu-a}, \cite{Abu-b}, \cite{KN2011},
\cite{Gol1999}).

\begin{proposition}
\label{S-prop}Let $S$ and $T$ be semirings.

\begin{enumerate}
\item $\mathbb{S}_{S}$ is a variety \emph{(}in the sense of Universal
Algebra, \emph{i.e.} a class of objects which is closed under homomorphic
images, direct products and subobjects\emph{)}.

\item $\mathbb{S}_{S}$ is complete \emph{(i.e.} has all small limits\emph{)}%
, equivalently $\mathbb{S}_{S}$ has equalizers, pullbacks and products. The
kernel of an $S$-linear map $f:M\longrightarrow N$ is
\begin{equation*}
\mathrm{Ker}(f):=\mathrm{Eq}(f,0)=\{m\in M\mid f(m)=0\}.
\end{equation*}

\item $\mathbb{S}_{S}$ is cocomplete \emph{(i.e.} has all small colimits%
\emph{)}, equivalently $\mathbb{S}_{S}$ has coequalizers \emph{(pushouts)}
and products. The cokernel of an $S$-linear map $f:M\longrightarrow N$ is $%
\mathrm{CoKer}(f):=\mathrm{Coeq}(f,0)=N/f(M).$

\item Every monomorphism in $\mathbb{S}_{S}$ is injective. A morphism in $%
\mathbb{S}_{S}$ is surjective if and only if it is a regular epimorphism.

\item $\mathbb{S}_{S}$ is a Barr-exact category \emph{\cite{Bar1971}} \emph{(%
}see also \emph{\cite{AHS2004})} with canonical factorization system given
by $(\mathbf{RegEpi},\mathbf{Mono})=(\mathbf{Surj},\mathbf{Inj}).$

\item $S_{S}$ is a regular generator in $\mathbb{S}_{S}.$

\item For all right $S$-semimodule $M$ and $N,$ we have a natural
isomorphism of Abelian monoids $\mathrm{Hom}_{S}(\mathfrak{c}(M),N)\simeq
\mathrm{Hom}_{S}(M,N),$ \emph{i.e.} the embeddings $\mathbb{CS}%
_{S}\hookrightarrow \mathbb{S}_{S}$ is right adjoint to $\mathfrak{c}(-);$
so, $\mathbb{CS}_{S}$ is a reflective subcategory of $\mathbb{S}_{S}.$

\item For every $(S,T)$-bisemimodule $_{S}M_{T},$ we have functors $M\otimes
_{T}-:$ $_{T}\mathbb{S}\longrightarrow $ $_{S}\mathbb{S}$ and $-\otimes
_{S}M:\mathbb{S}_{S}\longrightarrow \mathbb{S}_{T}.$ Moreover, we have
adjoint pairs of functors%
\begin{equation*}
(M\otimes _{T}-,\mathrm{Hom}_{S-}(M,-))\text{ and }(-\otimes _{S}M,\mathrm{%
Hom}_{-T}(M,-)),
\end{equation*}%
whence $M\otimes _{T}-$ and $-\otimes _{S}M$ preserves colimits.

\item $(_{S}\mathbb{S}_{S},\otimes _{S},S)$ is a biclosed monoidal category.

\item If $S$ is commutative, then $(\mathbb{S}_{S},\otimes _{S},S;\mathbf{%
\tau })$ is a symmetric braided monoidal category where $\mathbf{\tau }$ is
the \emph{flipping natural transformation}%
\begin{equation*}
M\otimes _{S}N\overset{\mathbf{\tau }_{M,N}}{\simeq }N\otimes _{S}M,\text{ }%
m\otimes _{S}n\longmapsto n\otimes _{S}m.
\end{equation*}
\end{enumerate}
\end{proposition}

\begin{definition}
Let $M$ and $N$ be $S$-semimodules. We call an $S$-linear map $%
f:M\longrightarrow N:$

$i$\emph{-uniform} (\emph{image-uniform}) iff $f(M)=\overline{f(M)};$

$k$\emph{-uniform} (\emph{kernel-uniform}) iff for all $m,m^{\prime }\in M$
we have
\begin{equation}
f(m)=f(m^{\prime })\Rightarrow m+k=m^{\prime }+k^{\prime }\text{ for some }%
k,k^{\prime }\in \mathrm{Ker}(f);  \label{k-uniform}
\end{equation}

\emph{uniform} iff $f$ is $i$-uniform and $k$-uniform.

We call $L\leq _{S}M$ a \emph{uniform subsemimodule} iff the embedding $L%
\overset{\iota _{L}}{\hookrightarrow }M$ is ($i$-)uniform, or equivalently
iff $L\leq _{S}M$ is subtractive. If $\equiv $ is an $S$-congruence on $M$
\cite{Gol1999}, then we call $M/\equiv $ a \emph{uniform quotient} iff the
projection $\pi _{\equiv }:M\longrightarrow M/\equiv $ is ($k$-)uniform.
\end{definition}

\begin{punto}
(\cite{Abu-b}) We say that a sequence $X\overset{f}{\longrightarrow }Y%
\overset{g}{\longrightarrow }Z$ of $S$-semimodules is \emph{exact} (resp.
\emph{semi-exact}, \emph{proper-exact, quasi-exact}) iff $f(X)=\mathrm{Ker}%
(g)$ and $g$ is $k$-uniform (resp. $\overline{f(X)}=\mathrm{Ker}(g),$ $f(X)=%
\mathrm{Ker}(g),$ $\overline{f(X)}=\mathrm{Ker}(g)$ and $g$ is $k$-uniform).
A (possibly infinite) sequence of $S$-semimodules $\cdots \longrightarrow
X_{i-1}\overset{f_{i-1}}{\longrightarrow }X_{i}\overset{f_{i}}{%
\longrightarrow }X_{i+1}\overset{f_{i+1}}{\longrightarrow }%
X_{i+2}\longrightarrow \cdots $ is said to be exact (resp. semi-exact,
proper exact, quasi-exact) iff each three-term subsequence $X_{i-1}\overset{%
f_{i-1}}{\longrightarrow }X_{i}\overset{f_{i}}{\longrightarrow }X_{i+1}$ is
exact (resp. semi-exact, proper exact, quasi-exact).
\end{punto}

\begin{lemma}
\label{exact}\emph{(\cite[Lemma 2.7]{Abu-b}) }Let $X,Y$ and $Z$ be $S$%
-semimodules.

\begin{enumerate}
\item $0\longrightarrow X\overset{f}{\longrightarrow }Y$ is exact if and
only if $f$ is injective.

\item $Y\overset{g}{\longrightarrow }Z\longrightarrow 0$ is exact if and
only if $g$ is surjective.

\item $0\longrightarrow X\overset{f}{\longrightarrow }Y\overset{g}{%
\longrightarrow }Z$ is semi-exact and $f$ is uniform if and only if $f$
induces an isomorphism $X\simeq \mathrm{Ker}(g).$

\item $X\overset{f}{\longrightarrow }Y\overset{g}{\longrightarrow }%
Z\longrightarrow 0$ is semi-exact and $g$ is uniform if and only if $g$
induces an isomorphism $Z\simeq \mathrm{Coker}(f).$

\item $0\longrightarrow X\overset{f}{\longrightarrow }Y\overset{g}{%
\longrightarrow }Z\longrightarrow 0$ is exact if and only if $f$ induces an
isomorphism $X\simeq \mathrm{Ker}(g)$ and $g$ induces an isomorphism $%
Z\simeq \mathrm{Coker}(f).$
\end{enumerate}
\end{lemma}

\begin{lemma}
\label{1st-iso}

\begin{enumerate}
\item An $S$-linear map $f:M\longrightarrow N$ induces an isomorphism $M/%
\mathrm{Ker}(f)\simeq f(M)$ if and only if $f$ is $k$-uniform.

\item For every $L\leq _{S}M,$ we have an exact sequence of $S$-semimodules%
\begin{equation*}
0\longrightarrow \overline{L}\overset{\iota _{\overline{L}}}{\longrightarrow
}M\overset{\pi _{L}}{\longrightarrow }M/L\longrightarrow 0.
\end{equation*}
\end{enumerate}
\end{lemma}

Since $\mathbb{S}_{S}$ is $(\mathbf{Surj},\mathbf{Inj})$\emph{-structured}
\cite{AHS2004} (and not $(\mathbf{Epi},\mathbf{Mono})$-structured), the
natural notions of projective objects, generators \emph{etc.} in this
category are defined relative to the class $\mathbf{Surj}$ of surjective $S$%
-linear maps (= regular epimorphisms) rather than the class $\mathbf{Epi}$
of all epimorphisms.

\begin{definition}
We say that an $S$-semimodule $X$ (\emph{uniformly}) \emph{generates} $M_{S}$
iff there exists an index set $\Lambda $ and a (uniform) surjective $S$%
-linear map $X^{(\Lambda )}\overset{\pi }{\longrightarrow }M\longrightarrow
0.$ With $\mathrm{Gen}(X)$ we denote the class of $S$-semimodules generated
by $X_{S}.$
\end{definition}

\begin{definition}
We say that $M_{S}$ is \emph{uniformly\emph{\ }}(\emph{\emph{finitely}})%
\emph{\ generated} iff there exists a (finite) index set $\Lambda $ and a
uniform surjective $S$-linear map $S^{(\Lambda )}\longrightarrow
M\longrightarrow 0.$
\end{definition}

\begin{remark}
Every $S$-semimodule $M$ is generated by $S:$ there exists a surjective $S$%
-linear map $S^{(\Lambda )}\overset{\pi }{\longrightarrow }M\longrightarrow
0 $ \cite[Proposition 17.11]{Gol1999}. However, it is not guaranteed that we
can find $\Lambda $ for which $\pi $ is uniform. Uniformly generated
semimodules were called $k$-semimodules in \cite{Alt1996}; we prefer the
terminology introduced above since it is more informative. Takahashi \cite%
{Tak1983} defined an $S$-semimodule $X$ to be \emph{normal} iff there exists
a projective $S$-semimodule $P$ and a uniform surjective $S$-linear map $P%
\overset{g}{\longrightarrow }X\longrightarrow 0$ (called a \emph{projective
presentation of }$X$). Indeed, every uniformly generated $S$-semimodule is
normal.
\end{remark}

\begin{punto}
\label{sgm-M}Let $M$ be a right $S$-semimodule. With $\sigma \lbrack M_{S}]$
($\sigma _{u}[M_{S}]$) we denote the closure of $\mathrm{Gen}(M_{S})$ under (%
\emph{uniform}) $S$-subsemimodules, \emph{i.e.} the smallest full
subcategory of $\mathbb{S}_{S}$ which contains $M_{S}$ and is closed under
direct sums, homomorphic images and (uniform) $S$-subsemimodules. We say
that $M_{S}$ is a (uniformly) \emph{subgenerator} for $\sigma \lbrack M_{S}]$
($\sigma _{u}[M_{S}]$). Notice that $\mathrm{Gen}(M_{S})\subseteq \sigma
_{u}[M_{S}]\subseteq \sigma \lbrack M_{S}].$
\end{punto}

\begin{remark}
Let $M$ be an $S$-semimodule and notice that the uniform $S$-subsemimodules
are precisely the kernels of $S$-linear maps by Lemma \ref{1st-iso} (2). It
follows that $X\in \sigma _{u}[M_{S}]$ if and only if $X\simeq \mathrm{Ker}%
(g),$ where $g:Y\longrightarrow Z$ is $S$-linear and $Y,Z\in \mathrm{Gen}%
(M_{S}),$ or equivalently if and only if there exist exact sequences of $S$%
-semimodules $0\longrightarrow X\overset{f}{\longrightarrow }Y\overset{g}{%
\longrightarrow }Z\longrightarrow 0$ in which $Y,Z\in \mathrm{Gen}(M_{S}).$
\end{remark}

\begin{proposition}
\label{sg=all}\emph{(cf. \cite[15.4]{Wis1991})} Let $M$ be a faithful $S$%
-semimodule and $T=\mathrm{End}(M_{S}).$ If $_{T}M$ is finitely generated,
then $\sigma \lbrack M_{S}]=\mathbb{S}_{S}.$
\end{proposition}

\begin{definition}
Let $\mathfrak{A}$ be a category with finite limits and finite colimits. A
functor $F:\mathfrak{A}\longrightarrow \mathfrak{B}$ is said to be \emph{%
left-exact}\textbf{\ }(\emph{right-exact}\textbf{)} iff $F$ preserves finite
limits (finite colimits). Moreover, $F$ is called \emph{exact} iff $F$ is
left-exact and right-exact.
\end{definition}

\begin{remarks}
Let $M$ be a right $S$-semimodules.

\begin{enumerate}
\item The \emph{contravariant functor} $\mathrm{Hom}_{S}(-,M):\mathbb{S}%
_{S}\longrightarrow \mathbf{AbMonoid}$ is left exact, whence it converts all
finite colimits into finite limits, \emph{e.g.} it converts coequalizers
into equalizers and converts pushouts into pullbacks. In particular, it
sends cokernels to kernels; this explains \cite[Theorem 2.6 (2)]{Tak1982} in
light of Lemma \ref{exact}.

\item The covariant functor $\mathrm{Hom}_{S}(M,-):\mathbb{S}%
_{S}\longrightarrow \mathbf{AbMonoid}$ is left exact, whence it preserves
all finite limits (\emph{e.g.} it sends equalizers to equalizers and
pullbacks to pullbacks). In particular, it preserves kernels; this explains
\cite[Theorem 2.6 (1)]{Tak1982} in light of Lemma \ref{exact}.

\item The covariant functor $M\otimes _{S}-:$ $_{S}\mathbb{S}\longrightarrow
\mathbf{AbMonoid}$ is right exact (since it has a right adjoint) whence it
preserves all finite colimits, \emph{e.g.} it sends coequalizers to
coequalizers and pushouts to pushouts. In particular, it sends cokernels to
kernels; this explains the \emph{analog} of \cite[Theorem 5.5]{Tak1982} in
light of Lemma \ref{exact}.
\end{enumerate}
\end{remarks}

\begin{definition}
We say that $M_{S}$ is

\emph{injective} iff for every monomorphism of $S$-semimodules (i.e.
injective $S$-linear map) $X\overset{f}{\longrightarrow }Y,$ every $S$%
-linear map $h:X\longrightarrow M$ can be extended to an $S$-linear map $%
\widetilde{h}:Y\longrightarrow M$ (such that $\widetilde{h}\circ f=h$);

\emph{uniformly injective} iff $\mathrm{Hom}_{S}(-,M):\mathbb{S}%
_{S}\longrightarrow \mathbf{AbMonoid}$ converts uniform monomorphisms into
uniform surjective maps (equivalently, $\mathrm{Hom}_{S}(-,M)$ preserves
short exact sequences);

$u$\emph{-injective} iff $\mathrm{Hom}_{S}(-,M):\mathbb{S}%
_{S}\longrightarrow \mathbf{AbMonoid}$ sends (uniform) monomorphisms to
(uniform) surjective maps.
\end{definition}

\begin{definition}
We say that $M_{S}$ is

\emph{projective} iff for every \emph{surjective} $S$-linear map $Y\overset{g%
}{\longrightarrow }Z\longrightarrow 0$ and every $S$-linear map $%
h:M\longrightarrow Z,$ there exists an $S$-linear map $\widetilde{h}%
:M\longrightarrow Y$ such that $h=g\circ \widetilde{h};$

\emph{uniformly projective} iff $\mathrm{Hom}_{S}(M,-):\mathbb{S}%
_{S}\longrightarrow \mathbf{AbMonoid}$ preserves uniform surjective
morphisms (equivalently, iff it preserves short exact sequences);

$u$\emph{-projective} iff $\mathrm{Hom}_{S}(M,-)$ sends (uniform) surjective
morphisms to (uniform) surjective maps.
\end{definition}

\begin{definition}
We say that $M_{S}$ is

\emph{cogenerator} iff for every $N_{S},$ there exist an index set $\Lambda $
and an $S$-linear embedding $N\hookrightarrow M^{\Lambda };$

\emph{uniformly} \emph{cogenerator} iff $\mathrm{Hom}_{S}(-,M):\mathbb{S}%
_{S}\longrightarrow \mathbf{AbMonoid}$ reflects short exact sequences.
\end{definition}

\begin{lemma}
\label{basis}

\begin{enumerate}
\item The following are equivalent for an $S$-semimodule $P_{S}:$

\begin{enumerate}
\item $P_{S}$ is projective;

\item $P_{S}$ is a \emph{retract} of a free $S$-semimodule, \emph{i.e.}
there exists an index set $\Lambda ,$ a surjective $S$-linear maps $%
S^{(\Lambda )}\overset{g}{\underset{f}{\rightleftarrows }}P$ with $g\circ f=%
\mathrm{id}_{P}.$

\item $P_{S}$ has a dual basis: there exists a subset $\{(p_{\lambda
},f_{\lambda })\}\subseteq P\times P^{\ast }$ such that:

\begin{itemize}
\item[-] For each $p\in P,$ the set $\Lambda (p)=\{\lambda \mid f_{\lambda
}(p)\neq 0\}$ is finite.

\item[-] $p=\sum p_{\lambda }f_{\lambda }(p).$
\end{itemize}
\end{enumerate}

\item If $P_{S}$ is uniformly generated and uniformly projective, then $P_{S}
$ is projective.
\end{enumerate}
\end{lemma}

\begin{Beweis}

\begin{enumerate}
\item $(a)\Longleftrightarrow (b)$ This is \cite[Proposition 17.16]{Gol1999}.

$(c)\Longleftrightarrow (d)$ The proof is similar to that of \cite[18.6]%
{Wis1991}, it can be shown that every projective $S$-semimodule has a dual
basis, whence finitely projective.

\item Let $S^{(\Lambda )}\overset{\pi }{\longrightarrow }P\longrightarrow 0$
be a uniform presentation of $P_{S}.$ Considering $\mathrm{id}%
_{P}:P\longrightarrow P,$ we find an $S$-linear map $h:P\longrightarrow
S^{(\Lambda )}$ such that $\pi \circ h=\mathrm{id}_{P},$ \emph{i.e.} $P_{S}$
is a retract of a free $S$-semimodule and so $P_{S}$ is projective by (1).$%
\blacksquare $
\end{enumerate}
\end{Beweis}

\begin{definition}
Let $M_{S}$ be an $S$-semimodule and set $M^{\ast }:=\mathrm{Hom}_{S}(M,S).$
We say that $M_{S}$ is \emph{finitely projective} iff for every finite
subset $\{m_{1},\cdots ,m_{l}\}\subseteq M,$ there exists $%
\{(m_{i},f_{j})\}_{j=1}^{n}\subset M\times M^{\ast }$ such that $%
m_{i}=\sum_{j=1}^{n}m_{j}f_{j}(m_{i})$ for each $i=1,\cdots ,l;$
\end{definition}

\begin{definition}
(\cite{Kat2004}, \cite{Abu-a}) We call a right $S$-semimodule $M:$

\emph{flat} iff $M\otimes _{A}-$ is left exact, \emph{i.e.} it preserves
finite limits, equivalently $M\simeq \lim\limits_{\longrightarrow
}F_{\lambda },$ a filtered limit of finitely generated free right $S$%
-semimodules;

\emph{uniformly flat} iff $M\otimes _{A}-:$ $_{A}\mathbb{S}\longrightarrow
\mathbf{AbMonoid}$ preserves \emph{uniform} subobjects;

\emph{mono-flat} iff $M\otimes _{A}-:$ $_{A}\mathbb{S}\longrightarrow
\mathbf{AbMonoid}$ preserves monomorphisms (injective $S$-linear maps);

$u$\emph{-flat} iff $M\otimes _{A}-:$ $_{A}\mathbb{S}\longrightarrow \mathbf{%
AbMonoid}$ sends (uniform) monomorphisms to (uniform) monomorphisms.
\end{definition}

\begin{remark}
\label{kernel-flat}Let $M$ be a right $S$-semimodule. Since $M\otimes _{S}-:$
$_{S}\mathbb{S}\longrightarrow \mathbf{AbMonoid}$ preserves direct sums, we
conclude that $M_{S}$ is flat if and only if $M\otimes _{S}-$ preserves
equalizers. Moreover, $M_{S}$ is flat if and only if $M\otimes _{S}-$
preserves pullbacks \cite{Kat2004}; so, flat semimodules are mono-flat. On
the other hand, if $M_{S}$ is a mono-flat $S$-semimodule, then $M_{S}$ is
uniformly flat (whence $u$-flat) if and only if $M\otimes _{S}-:$ $_{S}%
\mathbb{S}\longrightarrow \mathbf{AbMonoid}$ preserves kernels. It is not
known yet if there are examples of uniformly flat semimodules which are not
flat.
\end{remark}

\begin{definition}
Let $M$ be a right (left) $S$-semimodule. We say that a (uniform) $S$%
-subsemimodule $L\leq _{S}M$ is (\emph{uniformly})\emph{\ }$W$\emph{-pure}
for some left (right) $S$-semimodule $W$ iff $L\otimes _{S}W\leq M\otimes
_{S}W$ is a (uniform) submonoid. We call $L\leq _{S}M$ (\emph{uniformly})%
\emph{\ pure }iff $L\leq _{S}M$ is (uniformly) $W$-pure for every left
(right) $S$-semimodule $W.$ If $M$ is an $(S,T)$-bisemimodule and $L\leq
_{(S,T)}M,$ then we call $L\leq $ $_{(S,T)}M$ (\emph{uniformly})\emph{\ pure}
iff $L\hookrightarrow M$ is (uniformly) pure as a left $S$-subsemimodule as
well as a right $T$-subsemimodule.
\end{definition}

The proof of the following result is along the lines of \cite[Proposition 3.6%
]{Bou1974}:

\begin{lemma}
\label{Bou}Let $L$ be a right $S$-semimodule, $N$ a left $S$-semimodule, $%
K\leq _{S}L,$ $M\leq _{S}N$ and consider the exact sequences of $S$%
-semimodules%
\begin{equation*}
0\longrightarrow \overline{K}\overset{\iota _{\overline{K}}}{\longrightarrow
}L\overset{\pi _{K}}{\longrightarrow }L/K\longrightarrow 0\text{ and }%
0\longrightarrow \overline{M}\overset{\iota _{\overline{M}}}{\longrightarrow
}N\overset{\pi _{M}}{\longrightarrow }N/M\longrightarrow 0.
\end{equation*}%
We have an exact sequence of Abelian monoids%
\begin{equation*}
0\longrightarrow \overline{(\iota _{\overline{K}}\otimes _{S}N)(\overline{K}%
\otimes _{S}N)+(L\otimes _{S}\iota _{\overline{M}})(L\otimes _{S}\overline{M}%
)}\overset{\iota }{\longrightarrow }L\otimes _{S}N\overset{\pi _{K}\otimes
_{S}\pi _{M}}{\longrightarrow }L/K\otimes _{S}N/M\longrightarrow 0.
\end{equation*}%
If, moreover, $\overline{K}\leq _{S}L$ is $N$-pure and $\overline{M}\leq
_{S}N$ is $L$-pure, then we have an exact sequence of Abelian monoids%
\begin{equation*}
0\longrightarrow \overline{\overline{K}\otimes _{S}N+L\otimes _{S}\overline{M%
}}\overset{\iota }{\longrightarrow }L\otimes _{S}N\overset{\pi _{K}\otimes
_{S}\pi _{M}}{\longrightarrow }L/K\otimes _{S}N/M\longrightarrow 0.
\end{equation*}
\end{lemma}

\begin{definition}
We say that an $S$-semimodule $X$ is

\emph{finitely presented} iff $\mathrm{Hom}_{S}(X,-):S_{S}\longrightarrow
\mathbf{AbMonoid}$ preserves directed colimits (\emph{i.e.} $X\in \mathbb{S}%
_{S}$ is a finitely presentable object in the sense of \cite{AP1994});

\emph{uniformly finitely presented} iff $X$ is uniformly finitely generated
and for any exact sequence of $S$-semimodules%
\begin{equation*}
0\longrightarrow K\overset{f}{\longrightarrow }S^{n}\overset{g}{%
\longrightarrow }X\longrightarrow 0,
\end{equation*}%
the $S$-semimodule $K$ ($\simeq \mathrm{Ker}(g)$) is finitely generated.
\end{definition}

\begin{notation}
Let $M$ be a right $S$-semimodule. For every family $F=\{X_{\lambda
}\}_{\Lambda }$ of left $S$-semimodules, we have a morphism of Abelian
monoids%
\begin{equation*}
\varphi _{(M;F)}:M\otimes _{S}\dprod\limits_{\lambda \in \Lambda }X_{\lambda
}\longrightarrow \dprod\limits_{\lambda \in \Lambda }(M\otimes
_{S}X_{\lambda }),\text{ }m\otimes _{S}\{f_{\lambda }\}_{\Lambda }\mapsto
\{m\otimes _{S}f_{\lambda }\}_{\Lambda }.
\end{equation*}%
If $X_{\lambda }=S=X_{\gamma }$ for all $\lambda ,\gamma \in \Lambda ,$ then
we set $\widetilde{\varphi }_{M}=\varphi _{(M;S^{\Lambda })}:M\otimes
_{S}X^{\Lambda }\longrightarrow (M\otimes _{S}X)^{\Lambda }.$
\end{notation}

In the following lemma, we collect some properties of finitely presented
semimodules over semirings.

\begin{lemma}
\label{f.p.}Let $M$ be a right $S$-semimodule.

\begin{enumerate}
\item If $\widetilde{\varphi }_{M}:M\otimes _{S}S^{\Lambda }\longrightarrow
M^{\Lambda }$ is surjective for every index set $\Lambda ,$ then $M_{S}$ is
finitely generated.

\item If $\widetilde{\varphi }_{M}:M\otimes _{S}S^{\Lambda }$ is bijective
and $-\otimes _{S}S^{\Lambda }:\mathbb{S}_{S}\longrightarrow \mathbf{AbMonoid%
}$ preserves $i$-uniform morphisms for every index set $\Lambda ,$ then $%
M_{S}$ is uniformly finitely presented.

\item If $M_{S}$ is uniformly finitely presented, then $M_{S}$ has a finite
presentation through an exact sequence of $S$-semimodules%
\begin{equation}
S^{m}\overset{f}{\longrightarrow }S^{n}\overset{g}{\longrightarrow }%
M\longrightarrow 0.  \label{f.prese}
\end{equation}

\item If $M_{S}$ is uniformly finitely presented, then $M_{S}$ is finitely
presented.

\item If $M_{S}$ is finitely presented and flat, then $M_{S}$ is projective.
\end{enumerate}
\end{lemma}

\begin{Beweis}
The proofs are similar to the proofs for modules over rings.

\begin{enumerate}
\item The proof is similar to that of \cite[12.9 (1)]{Wis1991}.

\item The proof is similar to that of \cite[12.9 (2)]{Wis1991} using a
restricted version of the Short Five Lemma \cite[Lemma 1.22]{Abu-a} for
semimodules over semirings.

\item This is \cite[Proposition 2.25]{Abu-a}.

\item Given an arbitrary directed system $\{X_{i},\{f_{ij}\}\}_{I}$ of $S$%
-semimodules, we apply the contravariant functor $\mathrm{Hom}_{S}(-,%
\underset{\longrightarrow }{\lim }X_{i})$ to a finite presentation (\ref%
{f.prese}) of $M_{S},$ and then use a restricted version of the Short Five
Lemma \cite[Lemma 1.22]{Abu-a} for semimodules over semirings to prove that $%
\mathrm{Hom}_{S}(M,\underset{\longrightarrow }{\lim }X_{i})\simeq $ $%
\underset{\longrightarrow }{\lim }\mathrm{Hom}_{S}(M,X_{i}).$

\item Assume that $M_{S}$ is a finitely presented flat $S$-semimodule. By
definition of flat semimodules, $M=$ $\underset{\longrightarrow }{\lim }%
M_{i},$ where $\{M_{i}\}_{I}$ is a directed system of \emph{free }(\emph{%
projective})\emph{\ }$S$-semimodules. Since $M_{S}$ is finitely presented,
we have $\mathrm{End}_{S}(M)=\mathrm{Hom}_{S}(M,\underset{\longrightarrow }{%
\lim }M_{i})=$ $\underset{\longrightarrow }{\lim }\mathrm{Hom}_{S}(M,M_{i})$
and so $\mathrm{id}_{M}$ factorizes through some $M_{i},$ whence a retract
of the projective $S$-semimodule $M_{i}$ (cf. \cite[Proof of Theorem 2.6]%
{BR2004}). Since a retract of a projective $S$-semimodule is projective, we
conclude that $M_{S}$ is projective.$\blacksquare $
\end{enumerate}
\end{Beweis}

\section{Semicorings}

\qquad In this section, we introduce and investigate semicorings over
semirings and their categories of semicomodules.

Throughout, $S$ is a commutative semiring with $1_{S}\neq 0_{S}$ and $(%
\mathbb{S}_{S},\otimes _{S},S)$ is the symmetric monoidal category of $S$%
-semimodules \cite{Abu-a}. Moreover, $A$ is an $S$-semialgebra, \emph{i.e.}
a \emph{monoid} in $\mathbb{S}_{S},$ or equivalently a semiring $A$ along
with a morphism of semirings $\eta _{A}:S\longrightarrow A.$ With $_{A}%
\mathbb{S}$ ($\mathbb{S}_{A}$), we denote the category of left (right) $A$%
-semimodules and with $(_{A}\mathbb{S}_{A},\otimes _{A},A)$ the monoidal
category of $(A,A)$-bisemimodules.

\begin{punto}
By{\normalsize \ }an $A$\emph{-semiring}\ we mean a monoid in $_{A}\mathbb{S}%
_{A},$ \emph{i.e.} triple $(\mathcal{A},\mu _{\mathcal{A}},\eta _{\mathcal{A}%
})$ in which $\mathcal{A}$ is an $(A,A)$-bisemimodule and $\mu \mathcal{_{A}}%
:\mathcal{A}\otimes _{A}\mathcal{A}\longrightarrow \mathcal{A},$ $\eta _{%
\mathcal{A}}:A\longrightarrow \mathcal{A}$ are $(A,A)$-bilinear maps such
that the following diagrams are commutative%
\begin{equation*}
\begin{array}{ccc}
\xymatrix{\mathcal{A} \otimes_A \mathcal{A} \otimes_A \mathcal{A}
\ar[rr]^{\mu_{\mathcal A} \otimes_A \mathcal{A}} \ar[dd]_{\mathcal{A}
\otimes_A \mu_{\mathcal A} } & & \mathcal{A} \otimes_A \mathcal{A}
\ar[dd]^{\mu_{\mathcal A}}\\ & & \\ \mathcal{A} \otimes_A \mathcal{A}
\ar[rr]_{\mu_{\mathcal A}} & & \mathcal{A}} &  & \xymatrix{A \otimes_A
\mathcal{A} \ar[ddrr]_{\vartheta _{\mathcal A}^{l}} \ar[rr]^{ \eta_{\mathcal
A} \otimes_A \mathcal{A} } & & \mathcal{A} \otimes_A \mathcal{A}
\ar[dd]_{\mu_{\mathcal A}} & & \mathcal{A} \otimes_A A \ar[ll]_{{\mathcal A}
\otimes_A \eta_{\mathcal A}} \ar[ddll]^{\vartheta_{\mathcal A}^r} \\ & & & &
& \\ & & \mathcal{A} & & }%
\end{array}%
\end{equation*}%
We call $\mu _{\mathcal{A}}$ the \emph{multiplication }and $\eta \mathcal{%
_{A}}$ the \emph{unity} of $\mathcal{A}.$ If $A$ is commutative and $%
\mathcal{A}$ is an $A$-semiring with $xa=ax$ for all $x\in \mathcal{A}$ and $%
a\in A,$ then $\mathcal{A}$ is an $A$\emph{-semialgebra}. For $A$-semirings $%
\mathcal{A}$ and $\mathcal{B},$ we call an $(A,A)$-bilinear map $f:\mathcal{A%
}\longrightarrow \mathcal{B}$ a \emph{morphism of }$A$\emph{-semirings} iff $%
f\circ \mu _{\mathcal{A}}=\mu _{\mathcal{B}}\circ (f\otimes _{A}f)$ and $%
f\circ \eta _{\mathcal{A}}=\eta _{\mathcal{B}};$ the set of morphisms of $A$%
-semirings form $\mathcal{A}$ to $\mathcal{B}$ is denoted by $\mathrm{SRng}%
_{A}(\mathcal{A},\mathcal{B}).$ The category of $A$-semirings will be
denoted by $\mathbf{SRng}_{A}.$
\end{punto}

Corings over (associative) algebras were introduced by M. Sweedler \cite%
{Swe1975} as algebraic structures that are dual to rings. This suggests that
we define semicorings over semialgebras as algebraic structures dual to
semirings:

\begin{punto}
An\emph{\ }$A$\emph{-semicoring} is a \emph{comonoid} in $_{A}\mathbb{S}%
_{A}, $ equivalently a triple $(\mathcal{C},\Delta _{\mathcal{C}%
},\varepsilon _{\mathcal{C}})$ in which $\mathcal{C}$ is an $(A,A)$%
-bisemimodule and $\Delta _{\mathcal{C}}:\mathcal{C}\longrightarrow \mathcal{%
C}\otimes _{A}\mathcal{C}, $ $\varepsilon _{\mathcal{C}}:\mathcal{C}%
\longrightarrow A$ are $(A,A)$-bilinear maps such that the following
diagrams are commutative:%
\begin{equation*}
\begin{tabular}{lll}
$\xymatrix{{\mathcal {C}} \ar^(.45){\Delta_{\mathcal C}}[rr]
\ar_(.45){\Delta_{\mathcal C}}[d] & & {\mathcal {C}} \otimes_{A} {\mathcal
{C}} \ar^(.45){{\mathcal C} \otimes_A \Delta_{\mathcal C}}[d]\\ {\mathcal
{C}} \otimes_{A} {\mathcal {C}} \ar_(.45){\Delta_{\mathcal C} \otimes_A
{\mathcal C}}[rr] & & {\mathcal {C}} \otimes_{A} {\mathcal {C}} \otimes_{A}
{\mathcal {C}} }$ &  & $\xymatrix{ & {\mathcal {C}} \ar^(.45){\Delta
_{\mathcal {C}}}[d] & \\ {A} \otimes_{A} {\mathcal {C}}
\ar[ur]^(.45){\vartheta _{\mathcal {C}} ^l} & {\mathcal {C}} \otimes_{A}
{\mathcal {C}} \ar^(.45){\varepsilon _{\mathcal {C}} \otimes_A {\mathcal
C}}[l] \ar_(.45){{\mathcal C} \otimes_A \varepsilon _{\mathcal {C}}}[r] &
{\mathcal {C}} \otimes_{A} {A} \ar[ul]_(.45){\vartheta _{\mathcal {C}} ^r} }$%
\end{tabular}%
\end{equation*}%
We call $\Delta _{\mathcal{C}}$ the \emph{comultiplication} of $\mathcal{C}$
and $\varepsilon _{\mathcal{C}}$ the \emph{counity} of $\mathcal{C}.$
\end{punto}

\begin{punto}
For $A$-semicorings $(\mathcal{C},\Delta _{\mathcal{C}},\varepsilon _{%
\mathcal{C}}),$ $(\mathcal{D},\Delta _{\mathcal{D}},\varepsilon _{\mathcal{D}%
}),$ we call an $(A,A)$-bilinear map $f:\mathcal{D}\longrightarrow \mathcal{C%
}$\ an $A$\emph{-semicoring morphism} iff the following diagrams are
commutative%
\begin{equation*}
\begin{array}{ccc}
\xymatrix{{\mathcal {D}} \ar^(.45){f}[rr] \ar_(.45){\Delta_{\mathcal
{D}}}[d] & & {\mathcal {C}} \ar^(.45){\Delta_{{\mathcal {C}}}}[d] \\
{\mathcal {D}} \otimes_{A} {\mathcal {D}} \ar_(.45){f \otimes_A f}[rr] & &
{\mathcal {C}} \otimes_{A} {\mathcal {C}}} &  & \xymatrix{{\mathcal D}
\ar[dr]_{\varepsilon_{\mathcal D}} \ar[rr]^{f} & & {\mathcal C}
\ar[dl]^{\varepsilon_{\mathcal C}} \\ & A &}%
\end{array}%
\end{equation*}%
The set of $A$-semicoring morphisms from $\mathcal{D}$ to $\mathcal{C}$ is
denoted by $\mathrm{SCog}_{A}(\mathcal{D},\mathcal{C}).$ The category of $A$%
-semicorings is denoted by $\mathbf{SCorng}_{A}.$
\end{punto}

\begin{notation}
Let $(\mathcal{C},\Delta )$ be an $A$-semicoring. We use Sweedler-Heyneman's
$\sum $-notation, and write for $c\in \mathcal{C}:$%
\begin{equation*}
\Delta (c)=\sum c_{1}\otimes _{A}c_{2}\in \mathcal{C}\otimes _{A}\mathcal{C}.
\end{equation*}
\end{notation}

\begin{ex}
Let $\varphi :B\longrightarrow A$ be a morphism of $S$-semialgebras and
consider $A$ as a $(B,B)$-bisemimodule with actions given by $%
b\rightharpoonup a=\varphi (b)a$ and $a\leftharpoonup b=a\varphi (b)$ for
all $b\in B$ and $a\in A.$ Then $(\mathcal{C}:=A\otimes _{B}A,\Delta
,\varepsilon )$ is an $A$-semicoring where%
\begin{eqnarray*}
\Delta &:&\mathcal{C}\longrightarrow \mathcal{C}\otimes _{A}\mathcal{C},%
\text{ }a\otimes _{B}a^{\prime }\mapsto (a\otimes _{B}1_{A})\otimes
_{A}(1_{A}\otimes _{B}a^{\prime })=a\otimes _{B}1_{A}\otimes _{B}a^{\prime };
\\
\varepsilon &:&\mathcal{C}\longrightarrow A,\text{ }a\otimes _{B}a^{\prime
}\mapsto aa^{\prime }.
\end{eqnarray*}%
We call $(A\otimes _{B}A,\Delta ,\varepsilon )$ \emph{Sweedler's semicoring}.
\end{ex}

\begin{ex}
Let $M$ be an $(A,A)$-bisemimodule. We have an $A$-semicoring structure $%
\mathcal{C}=(A\oplus M,\Delta ,\varepsilon ),$ where%
\begin{eqnarray*}
\Delta &:&(a,m)\mapsto (a,0)\otimes _{A}(1,0)+(1,0)\otimes
_{A}(0,m)+(0,m)\otimes _{A}(1,0); \\
\varepsilon &:&(a,m)\mapsto a.
\end{eqnarray*}%
Notice that there are many properties $\mathbb{P}$ such that $_{A}\mathcal{C}
$ has Property $\mathbb{P}$ if (and only if) $_{A}M$ has Property $\mathbb{P}%
,$ \emph{e.g.} being flat, (finitely) projective, finitely generated \cite[%
Example 10 (1)]{Wisch1975}.
\end{ex}

\begin{punto}
Let $(C,\Delta ,\varepsilon )$ be a $S$-semicoring with $cs=sc$ for all $%
s\in S$ and $c\in C.$ We call $C$ an $S$\emph{-semicoalgebra}. An $S$%
-semicoalgebra is a \emph{comonoid} in the symmetric monoidal category $(%
\mathbb{S}_{S},\otimes _{S},S)$ of $S$-semimodules. If, moreover, $\sum
c_{1}\otimes _{S}c_{2}=\sum c_{2}\otimes _{S}c_{1}$ for all $c\in C,$ then
we say that $C$ is a \emph{cocommutative} $S$-semicoalgebra. We denote the
category of $S$-semicoalgebras by $\mathbf{SCoalg}_{S}$ and its \emph{full}
subcategory of cocommutative $S$-semicoalgebras by $_{\mathrm{coc}}\mathbf{%
SCoalg}_{S}.$
\end{punto}

\begin{ex}
Let $X$ be any set and consider the free $S$-semimodule with basis $X.$ Then
$(S^{(X)},\Delta ,\varepsilon )$ is an $S$-semicoalgebra where $\Delta $ and
$\varepsilon $ are defined by extending linearly the following assignments%
\begin{eqnarray*}
\Delta &:&S^{(X)}\mapsto S^{(X)}\otimes _{S}S^{(X)},\text{ }x\mapsto
x\otimes _{S}x; \\
\varepsilon &:&S^{(X)}\mapsto S,\text{ }x\mapsto 1_{S}.
\end{eqnarray*}%
Notice we have $S[X]\simeq S^{(\mathbb{N}_{0})}\simeq S[x],$ the polynomial
semiring in one indeterminate. So, $(S[x],\Delta _{1},\varepsilon _{1})$ is
an $S$-semicoalgebra with%
\begin{eqnarray*}
\Delta _{1} &:&S[x]\longrightarrow S[x]\otimes _{S}S[x],\text{ }%
\sum_{i=0}^{n}s_{i}x^{i}\mapsto \sum_{i=0}^{n}s_{i}x^{i}\otimes _{S}x^{i}; \\
\varepsilon _{1} &:&S[x]\longrightarrow S,\text{ }\sum_{i=0}^{n}s_{i}x^{i}%
\mapsto \sum_{i=0}^{n}s_{i}.
\end{eqnarray*}%
Moreover, $(S[x],\Delta _{2},\varepsilon _{2})$ is an $S$-semicoalgebra where%
\begin{eqnarray*}
\Delta _{2} &:&S[x]\longrightarrow S[x]\otimes _{S}S[x],\text{ }%
\sum_{i=0}^{n}s_{i}x^{i}\mapsto \sum_{i=0}^{n}s_{i}\left( \sum_{j=0}^{i}%
\binom{i}{j}x^{j}\otimes _{S}x^{j-i}\right) ; \\
\varepsilon _{2} &:&S[x]\longrightarrow S,\text{ }\sum_{i=0}^{n}s_{i}x^{i}%
\mapsto s_{0}.
\end{eqnarray*}
\end{ex}

\begin{ex}
(\cite{Wor2012}) Consider the idempotent Boolean semiring $\mathbf{B}%
=\{0,1\} $ (with $1+1=1\neq 0$). Let $P=\mathbf{B}<x,y\mid xy\neq yx>,$ the $%
\mathbf{B}$-semimodule of formal sums of words formed from the non-commuting
letters $x$ and $y.$ Then $P$ is in fact a $\mathbf{B}$-semialgebra with
multiplication given by the concatenation of words. It can be seen that $P$
has three structures of a $\mathbf{B}$-semicoalgebra given as follows:

\begin{enumerate}
\item $(P,\Delta _{1},\varepsilon _{1})$ is a $\mathbf{B}$-semicoalgebra
with $\Delta _{1}$ and $\varepsilon _{1}$ are defined on monomials and
extending linearly%
\begin{eqnarray*}
\Delta _{1} &:&P\longrightarrow P\otimes _{\mathbf{B}}P,\text{ }w\mapsto
w\otimes _{\mathbf{B}}w; \\
\varepsilon _{1} &:&P\longrightarrow \mathbf{B},\text{ }w\mapsto w(1,1).
\end{eqnarray*}

\item $(P,\Delta _{2},\varepsilon _{2})$ is a $\mathbf{B}$-semicoalgebra
with $\Delta _{2}$ and $\varepsilon _{2}$ are defined on monomials and
extended linearly%
\begin{eqnarray*}
\Delta _{2} &:&P\longrightarrow P\otimes _{\mathbf{B}}P,\text{ }w\mapsto
\sum_{w_{1}w_{2}=w}w_{1}\otimes _{\mathbf{B}}w_{2}; \\
\varepsilon _{2} &:&P\longrightarrow \mathbf{B},\text{ }w\mapsto w(0,0).
\end{eqnarray*}

\item $(P,\Delta _{3},\varepsilon _{3})$ is a $\mathbf{B}$-semicoalgebra
with $\Delta _{3}$ and $\varepsilon _{3}$ are defined on monomials and
extended as semialgebra morphisms%
\begin{eqnarray*}
\Delta _{3} &:&P\longrightarrow P\otimes _{\mathbf{B}}P,\text{ }\Delta
(x)=1\otimes _{\mathbf{B}}x+x\otimes _{\mathbf{B}}1,\text{ }\Delta
(y)=1\otimes _{\mathbf{B}}y+y\otimes _{\mathbf{B}}1; \\
\varepsilon _{3} &:&P\longrightarrow \mathbf{B},\text{ }w\mapsto w(0,0).
\end{eqnarray*}
\end{enumerate}
\end{ex}

In what follows, we mean by \emph{locally presentable} categories those in
the sense of \cite{AP1994}.

\begin{definition}
(\cite[Definition 1.17]{AP1994})\ Let $\mathfrak{A}$ be a category and $%
\lambda $ a regular cardinal. We say that an object $X\in \mathfrak{A}$ is
\emph{locally }$\lambda $\emph{-presentable} iff $\mathfrak{A}(X,-)$
preserves $\lambda $-directed colimits. A category $\mathfrak{A}$ is said to
be \emph{locally presentable} iff $\mathfrak{A}$ is cocomplete and has a set
$\mathbf{P}$ of $\lambda $-presentable objects, for some regular cardinal $%
\lambda ,$ such that every object in $\mathfrak{A}$ is a $\lambda $-directed
colimit of objects from $\mathbf{P}.$
\end{definition}

\begin{punto}
Applying results of Porst \cite{Por2008b} iteratively, we have the
following: The category $(\mathbf{Set},\times ,\{\ast \})$ of sets is an
\emph{admissible} symmetric monoidal category which is locally presentable.
It follows that the category $\mathbf{Monoid}=\mathbf{Mon}(\mathbf{Set})$ of
monoids is finitary monadic over $\mathbf{Set}$ and is locally presentable;
the full subcategory $\mathbf{AbMonoid}=$ $_{\mathrm{c}}\mathbf{Mon}(\mathbf{%
Set})$ of Abelian (commutative) monoids is finitary monadic over $\mathbf{Set%
}$ and locally presentable. Notice that $(\mathbf{AbMonoid},\otimes ,\mathbb{%
N}_{0})\simeq (\mathbb{S}_{\mathbb{N}_{0}},\otimes _{\mathbb{N}_{0}},\mathbb{%
N}_{0}),$ the category of semimodules over the semiring $\mathbb{N}_{0}$ of
non-negative integers, is a \emph{biclosed} symmetric monoidal category and
it follows that the category of semirings $\mathbf{SRng}=\mathbf{Mon}(%
\mathbf{AbMonoid})$ is finitary monadic over $\mathbf{AbMonoid}$ and locally
presentable; the full subcategory of \emph{commutative} semirings $_{\mathrm{%
c}}\mathbf{SRng}\simeq $ $_{\mathrm{c}}\mathbf{Mon}(\mathbf{AbMonoid})$ is
finitary monadic over $\mathbf{AbMonoid}$ and locally presentable. Moreover,
for any (commutative) semiring $A$ ($S$), the category $_{A}\mathbb{S}_{A}$ (%
$\mathbb{S}_{S}$) of $(A,A)$-bisemimodules ($S$-semimodules) is locally
presentable since it is a variety \cite[1.10 (2)]{AP1994}.
\end{punto}

The \emph{Fundamental Theorem of Coalgebras} \cite{Swe1969} states that
every coalgebra over a field is a directed limit of finite dimensional (%
\textrm{equivalently} locally presentable \cite[Proposition 1]{Por2006})
subcoalgebras. This result was generalized to comonoids in a locally
presentable symmetric monoidal category by Porst \cite{Por2008b} (see \cite%
{Por2006} and \cite{Por2008}). The results of Porst apply in particular to
semicoalgebras over commutative semirings.

\begin{proposition}
\label{semicoal}Consider the categories $\mathbf{SAlg}_{S}$ of $S$%
-semialgebras and $\mathbf{SCoalg}_{S}$ of $S$-semicoalgebras.

\begin{enumerate}
\item $\mathbf{SAlg}_{S}$ is finitary monadic over $\mathbb{S}_{S}$ and
locally presentable.

\item $_{\mathrm{c}}\mathbf{SAlg}_{S}$ is reflective in $\mathbf{SAlg}_{S},$
finitary monadic over $\mathbb{S}_{S}$ and locally presentable.

\item $_{\mathrm{c}}\mathbf{SAlg}_{S}$ is closed in $\mathbf{SAlg}_{S}$
under limits, directed colimits and absolute colimits\footnote{%
A limit (colimit) $\mathcal{K}$ is said to be \emph{absolute} iff $\mathcal{K%
}$ is preserved by every functor $G:\mathfrak{A}\longrightarrow \mathfrak{B}%
, $ where $\mathfrak{B}$ is an arbitrary category \cite[20.14 (3)]{AHS2004}.}%
.

\item $\mathbf{SCoalg}_{S}$ is comonadic over $\mathbb{S}_{S}$ and locally
presentable.

\item $_{\mathrm{coc}}\mathbf{SCoalg}_{S}$ is coreflective in $\mathbf{SCoalg%
}_{S},$ comonadic over $\mathbb{S}_{S}$ and locally presentable.

\item $_{\mathrm{coc}}\mathbf{SCoalg}_{S}$ is closed in $\mathbf{SCoalg}_{S}$
under colimits and absolute limits.
\end{enumerate}
\end{proposition}

\begin{Beweis}
The result is an immediate application of the main results of \cite{Por2008b}
taking into consideration that $(\mathbb{S}_{S},\otimes _{S},S)$ is a
biclosed (whence admissible) symmetric monoidal category and that we have
isomorphisms of categories $\mathbf{SAlg}_{S}\simeq \mathbf{Monoid}(\mathbb{S%
}_{S})$ and $\mathbf{SCoalg}_{S}\simeq \mathbf{Comonoid}(\mathbb{S}%
_{S}).\blacksquare $
\end{Beweis}

The proof of the following result is essentially the same as that for
corings over an algebra \cite{Por2006}.

\begin{proposition}
\label{coring}The category $\mathbf{SCoring}_{A}$ of $A$-semicorings is
comonadic over $_{A}\mathbb{S}_{A},$ locally presentable and a covariety
\emph{(}in the sense of \emph{\cite{AP2003})}.
\end{proposition}

\begin{punto}
Let $(\mathcal{C},\Delta _{\mathcal{C}},\varepsilon _{\mathcal{C}})$ be an $%
A $-semicoring and $\mathcal{D}\leq _{(A,A)}\mathcal{C}.$ We say that $%
\mathcal{D}$ is an $A$-\emph{subsemicoring} of $\mathcal{C}$ iff $\mathcal{D}
$ is an $A$-semicoring and the embedding $\iota :\mathcal{D}\hookrightarrow
\mathcal{C}$ is a morphism of $A$-semicorings. If $\mathcal{D}\leq _{(A,A)}%
\mathcal{C}$ is \emph{pure}, then $\mathcal{D}$ is an $A$\emph{-subsemicoring%
} of $\mathcal{C}$ if and only if $\Delta _{\mathcal{C}}(\mathcal{D}%
)\subseteq \mathcal{D}\otimes _{A}\mathcal{D}\subseteq \mathcal{C}\otimes
_{A}\mathcal{C}$ and $\varepsilon _{\mathcal{D}}$ is the restriction of $%
\varepsilon _{\mathcal{C}}$ to $\mathcal{D}.$
\end{punto}

\begin{punto}
Let $\mathcal{C}$ be a \emph{coassociative }$A$-semicoring. Associated to $%
\mathcal{C}$ are three dual $A$-semirings:

$^{\ast }\mathcal{C}:=(\mathrm{Hom}_{A-}(\mathcal{C},A),\star _{l})$ is an $%
A $-semiring, where%
\begin{equation*}
(f\star _{l}g)(c)=\sum g(c_{1}f(c_{2}))\text{ for all }f,g\in \text{ }^{\ast
}\mathcal{C}\text{ and }c\in \mathcal{C}{\normalsize ;}
\end{equation*}

$\mathcal{C}^{\ast }:=(\mathrm{Hom}_{-A}(\mathcal{C},A),\star _{r})$ is an $%
A $-semiring, where
\begin{equation*}
(f\star _{r}g)(c)=\sum f(g(c_{1})c_{2})\text{ for all }f,g\in \mathcal{C}%
^{\ast }\text{ and }c\in \mathcal{C}{\normalsize ;}
\end{equation*}

$^{\ast }\mathcal{C}^{\ast }:=(\mathrm{Hom}_{(A,A)}(\mathcal{C},A),\star )$
is an\emph{\ }$A$-semiring, where
\begin{equation*}
(f\star g)(c)=\sum g(c_{1})f(c_{2})\text{ for all }f,g\in \text{ }^{\ast }%
\mathcal{C}^{\ast }\text{ and }c\in \mathcal{C}{\normalsize .}
\end{equation*}%
The counity $\varepsilon _{\mathcal{C}}$ is a unity for $^{\ast }\mathcal{C}%
, $ $\mathcal{C}\mathbf{^{\ast }}$ and $^{\ast }\mathcal{C}^{\ast }.$
\end{punto}

\begin{definition}
Let $\mathcal{C}$ be an $A$-semicoring. We call

$K\leq \mathcal{C}_{A}$ a right $\mathcal{C}$-coideal iff $\Delta _{\mathcal{%
C}}(K)\subseteq \func{Im}(\iota _{K}\otimes _{A}\mathcal{C});$

$K\leq $ $_{A}\mathcal{C}$ a left $\mathcal{C}$-coideal iff $\Delta _{%
\mathcal{C}}(K)\subseteq \func{Im}(\mathcal{C}\otimes _{A}\iota _{K});$

$K\leq _{(A,A)}\mathcal{C}$ a $\mathcal{C}$-\emph{bicoideal} iff $\Delta _{%
\mathcal{C}}(K)\subseteq \func{Im}(\iota _{K}\otimes _{A}\mathcal{C})\cap
\func{Im}(\mathcal{C}\otimes _{A}\iota _{K});$

$K\leq _{(A,A)}\mathcal{C}$ a $\mathcal{C}$-\emph{coideal} iff $K=\mathrm{Ker%
}(f)$ for some \emph{uniform} surjective morphism of $A$-semicorings $f:%
\mathcal{C}\longrightarrow \mathcal{C}^{\prime }.$
\end{definition}

\begin{proposition}
\label{coideal}Let $(\mathcal{C},\Delta _{\mathcal{C}},\varepsilon _{%
\mathcal{C}})$ be an $A$-semicoring, $K\leq _{(A,A)}\mathcal{C}$ be uniform
and consider the canonical \emph{(}uniform\emph{)} surjection $\pi _{K}:%
\mathcal{C}\longrightarrow \mathcal{C}/K.$ The following are equivalent:

\begin{enumerate}
\item $K$ is a coideal of $\mathcal{C};$

\item There exists a morphism of $A$-semicorings $f:C\longrightarrow
C^{\prime }$ and an exact sequence of $(A,A)$-bisemimodules%
\begin{equation*}
0\longrightarrow K\overset{\iota }{\longrightarrow }\mathcal{C}\overset{f}{%
\longrightarrow }\mathcal{C}^{\prime }\longrightarrow 0;
\end{equation*}

\item $\mathcal{C}/K$ is an $A$-semicoring and the $\pi _{K}:\mathcal{C}%
\longrightarrow \mathcal{C}/K$ is a morphism of $A$-semicorings;

\item $\Delta _{\mathcal{C}}(K)\subseteq \overline{(\iota _{K}\otimes _{A}%
\mathcal{C})(K\otimes _{A}\mathcal{C})+(\mathcal{C}\otimes _{A}\iota _{K})(%
\mathcal{C}\otimes _{A}K)}$ and $\varepsilon _{\mathcal{C}}(K)=0.$

If $K\leq _{S}\mathcal{C}$ is uniformly $\mathcal{C}$-pure, then these are
moreover equivalent to:

\item[5] $\Delta _{\mathcal{C}}(K)\subseteq \overline{K\otimes _{A}\mathcal{C%
}+\mathcal{C}\otimes _{A}K}$ and $\varepsilon _{\mathcal{C}}(K)=0.$
\end{enumerate}
\end{proposition}

\begin{Beweis}
First of all, notice that the uniform subsemimodules are precisely the
subtractive ones, whence $\overline{K}=K.$

$(1)\Longleftrightarrow (2)$ Follows directly from the definition and Lemma %
\ref{exact}.

$(2)\Rightarrow (3)$ By Lemma \ref{1st-iso} (1), $f$ induces an isomorphism
of $(A,A)$-bisemimodules $\mathcal{C}/K\overset{\overline{f}}{\simeq }%
\mathcal{C}^{\prime }.$ One can easily check that this isomorphisms provides
$\mathcal{C}/K$ with a structure of an $A$-semicoring and that $\pi _{K}=%
\overline{f}^{-1}\circ f:\mathcal{C}\longrightarrow \mathcal{C}/K$ is a
morphism of $A$-semicorings.

$(3)\Rightarrow (2)$ Since $K\leq _{(A,A)}\mathcal{C}$ is uniform, it
follows by Lemma \ref{exact} that $K\simeq \mathrm{Ker}(\pi _{K}),$ whence $%
K $ is a coideal (notice that $\pi _{K}$ is uniform).

$(3)\Rightarrow (4)$ Consider the following diagram of $(A,A)$-bisemimodules%
\begin{equation}
\xymatrix{0 \ar[rr] & & K \ar[rr]^{\iota_K} \ar@{.>}[d]^{\kappa} & &
{\mathcal C} \ar[d]^{\Delta} \ar[rr]^{\pi_K} & & {\mathcal C}/K
\ar@{.>}[d]^{\overline{\Delta}} \ar[rr] & & 0\\ 0 \ar[rr] & & {\rm Ker}
(\pi_K \otimes_S \pi_K) \ar[rr]_{\iota} & & {\mathcal C} \otimes_S {\mathcal
C} \ar[rr]_{\pi_K \otimes_S \pi_K} & & {\mathcal C}/K \otimes_S {\mathcal
C}/K \ar[rr] & & 0}  \label{C/K}
\end{equation}%
By assumption, the second square is commutative and so $(\pi _{K}\otimes
_{A}\pi _{K})\circ \Delta _{\mathcal{C}}\circ \iota _{K}=\overline{\Delta }%
\circ \pi _{K}\circ \iota _{K}=0.$ By the universal property of kernels,
there exists an $(A,A)$-bilinear map $\kappa :K\longrightarrow \mathrm{Ker}%
(\pi _{K}\otimes _{A}\pi _{K})$ such that the first square is commutative,
equivalently $\Delta _{\mathcal{C}}(K)\subseteq \mathrm{Ker}(\pi _{K}\otimes
_{A}\pi _{K})=\overline{(\iota _{K}\otimes _{A}\mathcal{C})(K\otimes _{A}%
\mathcal{C})+(\mathcal{C}\otimes _{A}\iota _{K})(\mathcal{C}\otimes _{A}K)}$
by Lemma \ref{Bou}. Moreover, we have $\overline{\varepsilon }\circ \pi
_{K}=\varepsilon _{\mathcal{C}},$ whence $\varepsilon _{\mathcal{C}}(K)=(%
\overline{\varepsilon }\circ \pi _{K})(K)=0.$

$(4)\Rightarrow (3)$ Consider Diagram (\ref{C/K}). By assumption, the first
square is commutative and so $(\pi _{K}\otimes _{A}\pi _{K})\circ \Delta _{%
\mathcal{C}}\circ \iota _{K}=(\pi _{K}\otimes _{A}\pi _{K})\circ \iota \circ
\kappa =0.$ By the universal property of cokernels, there exists a unique $%
(A,A)$-bilinear map $\overline{\Delta }:C/K\longrightarrow C/K\otimes
_{A}C/K $ such that the second square is commutative. Moreover, since $%
\varepsilon _{\mathcal{C}}(K)=0,$ the assignment%
\begin{equation*}
\overline{\varepsilon }:\mathcal{C}/K\longrightarrow A,\overline{c}\mapsto
\varepsilon _{\mathcal{C}}(c)
\end{equation*}%
is a well defined $(A,A)$-bilinear map. One can easily check that $(\mathcal{%
C}/K,\overline{\Delta },\overline{\varepsilon })$ is an $A$-semicoring and
that $\pi _{K}:\mathcal{C}\longrightarrow \mathcal{C}/K$ is a morphism of $A$%
-semicorings.

If $K\leq _{(A,A)}\mathcal{C}$ is uniformly $\mathcal{C}$-pure, then $%
\mathrm{Ker}(\pi _{K}\otimes _{A}\pi _{K})=\overline{K\otimes _{A}\mathcal{C}%
+\mathcal{C}\otimes _{A}K}$ by Lemma \ref{Bou} and so the last assertion
follows.$\blacksquare $
\end{Beweis}

\subsection*{Semicomodules}

\qquad Dual to semimodules of semirings are semicomodules of semicorings:

\begin{punto}
Let $(\mathcal{C},\Delta ,\varepsilon )$ be an $A$-semicoring. A right $%
\mathcal{C}$\emph{-semicomodule} is a right $A$-semimodule $M$ associated
with an $A$-linear map (called $\mathcal{C}$\emph{-coaction})
\begin{equation*}
\rho ^{M}:M\longrightarrow M\otimes _{A}\mathcal{C},\text{ }m\mapsto \sum
m_{<0>}\otimes _{A}m_{<1>},
\end{equation*}%
such that the following diagrams are commutative%
\begin{equation*}
\begin{array}{ccc}
\xymatrix{M \ar^(.4){\rho ^M}[rr] \ar_(.45){\rho ^M}[d] & & M \otimes_{A}
{\mathcal{C}} \ar^(.45){M \otimes_A \Delta}[d] \\ M \otimes_{A}
{\mathcal{C}} \ar_(.4){\rho ^M \otimes_A {\mathcal{C}}}[rr] & & M
\otimes_{A} {\mathcal {C}} \otimes_{A} {\mathcal{C}}} &  & \xymatrix{M
\ar[rr]^{\rho^M} & & M \otimes_A {\mathcal C} \ar[dl]^{M \otimes_A
\varepsilon} \\ & M \otimes_A A \ar[ul]^{\vartheta_M^r} & }%
\end{array}%
\end{equation*}%
Let $M$ and $N${\normalsize \ }be right $\mathcal{C}$-semicomodules. We call
an $A$-linear map $f:M\longrightarrow N$ a $\mathcal{C}$\emph{-semicomodule
morphism}{\normalsize \ }(or $\mathcal{C}$\emph{-colinear}) iff the
following diagram is commutative
\begin{equation*}
\xymatrix{M \ar[rr]^{f} \ar[d]_{\rho ^M} & & N \ar[d]^{\rho ^N}\\ M
\otimes_{A} {\mathcal{C}} \ar[rr]_{f \otimes_A {\mathcal{C}}} & & N
\otimes_{A} {\mathcal{C}} }
\end{equation*}%
The set of $\mathcal{C}$-colinear maps from $M$ to $N$ is denoted by $%
\mathrm{Hom}^{\mathcal{C}}(M,N).$ The category of right $\mathcal{C}$%
-semicomodules and $\mathcal{C}$-colinear maps is denoted by $\mathbb{S}^{%
\mathcal{C}}.$ For a right $\mathcal{C}$-semicomodule $M,$ we call $L\leq
_{A}M$ a $\mathcal{C}$\emph{-subsemicomodule} iff $(L,\rho ^{L})\in \mathbb{S%
}^{\mathcal{C}}$ and the embedding $L\overset{\iota _{L}}{\hookrightarrow }M$
is $\mathcal{C}$-colinear. Symmetrically, we define the category $^{\mathcal{%
C}}\mathbb{S}$ of left $\mathcal{C}$-semicomodules. For two left $\mathcal{C}
$-semicomodules $M$ and $N,$ we denote by $^{\mathcal{C}}\mathrm{Hom}(M,N)$
the set of $\mathcal{C}$-colinear maps from $M$ to $N.$
\end{punto}

\begin{punto}
Let $(M,\rho ^{(M;\mathcal{C)}})$ be a right $\mathcal{C}$-semicomodule, $%
(M,\rho ^{(M;D\mathcal{)}})$ a left $\mathcal{D}$-semicomodule and consider
the left $\mathcal{D}$-semicomodule $(M\otimes _{A}\mathcal{C},\rho ^{(M;D%
\mathcal{)}}\otimes _{A}\mathcal{C})$ (the right $\mathcal{C}$-semicomodule $%
(\mathcal{D}\otimes _{A}M,\mathcal{D}\otimes _{A}\rho ^{(M;\mathcal{C)}})$).
We call $M$ a $(\mathcal{D},\mathcal{C})$\emph{-bisemicomodule} iff $\rho
^{(M;\mathcal{C)}}:M\longrightarrow M\otimes _{A}\mathcal{C}$ is $\mathcal{D}
$-colinear, or equivalently iff $\rho ^{(M;D\mathcal{)}}:M\longrightarrow
\mathcal{D}\otimes _{A}M$ is $\mathcal{C}$-colinear. For $(\mathcal{D},%
\mathcal{C})$-bisemicomodules $M$ and $N,$ we call a $\mathcal{D}$-colinear $%
\mathcal{C}$-colinear map $f:M\longrightarrow N$ a $(\mathcal{D},\mathcal{C}%
) $\emph{-bisemicomodule morphism}\textbf{\ }(or $(\mathcal{D},\mathcal{C})$%
\emph{-bicolinear}). The category of $(\mathcal{D},\mathcal{C})$%
-bisemicomodules and $(\mathcal{D},\mathcal{C})$-bicolinear maps is denoted
by $^{\mathcal{D}}\mathbb{S}^{\mathcal{C}}.$
\end{punto}

\begin{remark}
\label{inj-kounit}Let $(\mathcal{C},\Delta ,\varepsilon )$ be an $A$%
-semicoring. If $(M,\rho ^{M})$ is a right $\mathcal{C}$-semicomodule, then $%
\rho ^{M}$ is a splitting monomorphism in $\mathbb{S}_{A}$ (but $M$ is not
necessarily a direct summand of $M\otimes _{A}\mathcal{C}\mathbf{;}$ see
\cite[16.6]{Gol1999}).
\end{remark}

Although every $S$-semialgebra $A$ is a regular generator in $\mathbb{S}_{A}$
and in $_{A}\mathbb{S},$ it is not evident that $A$ is a generator in $_{A}%
\mathbb{S}_{A}$\ (even if $S$ is a commutative ring and $A$ is an $A$%
-algebra \cite[28.1]{Wis1996}.

\begin{definition}
We define the \emph{centroid }of the $S$-semialgebra $A$ as%
\begin{equation*}
C(A):=\{f\in \mathrm{End}_{S}(A)\mid af(b)=f(ab)=f(a)b\text{ for all }a,b\in
A\}.
\end{equation*}%
We say that $A$ is a \emph{central }$S$\emph{-semialgebra} iff $S\overset{%
\varphi }{\simeq }C(A),$ where%
\begin{equation*}
\varphi :S\longrightarrow C(A),\text{ }a\mapsto \lbrack a\mapsto sa].
\end{equation*}%
We say that $A$ is an \emph{Azumaya }$S$\emph{-semialgebra} iff $A$ is a
central $S$-semialgebra such that $A$ is a regular generator in $_{A}\mathbb{%
S}_{A}.$
\end{definition}

We present now the main reconstruction result:

\begin{theorem}
\label{comonad-C}

\begin{enumerate}
\item Let $\mathcal{C}$ be an $(A,A)$-bisemimodule. The following are
equivalent:

\begin{enumerate}
\item $\mathcal{C}$ is an $A$-semicoring;

\item $\mathcal{C}\otimes _{A}-:$ $_{A}\mathbb{S}\longrightarrow $ $_{A}%
\mathbb{S}$ is a comonad;

\item $-\otimes _{A}\mathcal{C}:\mathbb{S}_{A}\longrightarrow \mathbb{S}_{A}
$ is a comonad.
\end{enumerate}

\item If $A$ is an Azumaya $S$-semialgebra, then there is a bijective
correspondence between the structures of $A$-semirings on $\mathcal{C},$ the
comonad structures on $\mathcal{C}\otimes _{A}-:$ $_{A}\mathbb{S}%
\longrightarrow $ $_{A}\mathbb{S}$ and the comonad structures on $-\otimes
_{A}\mathcal{C}:$ $_{A}\mathbb{S}\longrightarrow $ $_{A}\mathbb{S}.$

\item Let $\mathcal{C}$ be an $A$-semicoring and $\mathcal{D}$ a $B$%
-semicoring \emph{(}for some $S$-semialgebra $B$\emph{)}. We have
isomorphisms of categories%
\begin{eqnarray*}
^{\mathcal{D}}\mathbb{S} &\simeq &\text{ }(_{B}\mathbb{S})^{_{\mathcal{D}%
\otimes _{B}-}},\text{ }\mathbb{S}^{\mathcal{C}}\simeq (\mathbb{S}%
_{A})^{-\otimes _{A}\mathcal{C}};\text{ } \\
((_{B}\mathbb{S}_{A})^{-\otimes _{A}\mathcal{C}})^{_{\mathcal{D}\otimes
_{B}-}} &\simeq &\text{ }^{\mathcal{D}}\mathbb{S}^{\mathcal{C}}\simeq \text{
}((_{B}\mathbb{S}_{A})^{_{\mathcal{D}\otimes _{B}-}})^{-\otimes _{A}\mathcal{%
C}}.
\end{eqnarray*}
\end{enumerate}
\end{theorem}

\begin{Beweis}
(1) and (3) follow directly from the definitions \cite[18.28]{BW2003}. The
proof of the bijective correspondence in (2) is similar to that of \cite[%
Theorem 3.9]{Ver} taking into consideration that $_{A}\mathbb{S}_{A}$ is
cocomplete, that $A$ is a regular generator in $_{A}\mathbb{S}_{A}$ (by our
assumption that $A$ is an Azumaya $S$-semialgebra) and the fact that $%
-\otimes _{A}X$ and $X\otimes _{A}-$ preserve colimits in $_{A}\mathbb{S}%
_{A} $ for every $(A,A)$-bisemimodule $X.\blacksquare $
\end{Beweis}

The main setting in \cite[p. 228]{Por2008} applies perfectly to our context.
In particular, the category $(_{A}\mathbb{S}_{A},\otimes _{A},A)$ is a
pointed monoidal category, $_{A}\mathbb{S}_{A}$ is a variety in the sense of
whence a locally presentable category \cite{AP1994}. Moreover, for each $%
_{A}X_{A},$ the functor $X\otimes _{A}-:$ $_{A}\mathbb{S}_{A}\longrightarrow
$ $_{A}\mathbb{S}_{A}$ has a right adjoint given by $\mathrm{Hom}_{A-}(X,-)$
and the functor $-\otimes _{A}X:$ $_{A}\mathbb{S}_{A}\longrightarrow $ $_{A}%
\mathbb{S}_{A}$ has a right adjoint given by $\mathrm{Hom}_{-A}(X,-).$ So,
the following result are essentially the same as in the proof of the
corresponding ones in \cite{Por2006}.

\begin{proposition}
\label{Porst}Let $\mathcal{C}$ be an $A$-semicoring and $\mathcal{F}:\mathbb{%
S}^{\mathcal{C}}\longrightarrow \mathbb{S}_{A}$ the forgetful functor.

\begin{enumerate}
\item $\mathbb{S}^{\mathcal{C}}$ is comonadic, locally presentable and a
covariety.

\item $\mathcal{F}$ creates all colimits and isomorphisms.

\item $\mathbb{S}^{\mathcal{C}}$ is cocomplete, \emph{i.e.} $\mathbb{S}^{%
\mathcal{C}}$ has all \emph{(}small\emph{)} colimits, \emph{e.g.}
coequalizers, cokernels, pushouts, directed colimits and direct sums.
Moreover, the colimits are formed in $\mathbb{S}_{A}.$

\item $\mathbb{S}^{\mathcal{C}}$ is complete, \emph{i.e.} $\mathbb{S}^{%
\mathcal{C}}$ has all \emph{(}small\emph{)} limits, \emph{e.g.} equalizers,
kernels, pullbacks, inverse limits and direct products. Moreover, $\mathcal{F%
}$ creates all limits preserved by $-\otimes _{A}\mathcal{C}:\mathbb{S}%
_{A}\longrightarrow \mathbb{S}_{A}.$
\end{enumerate}
\end{proposition}

\begin{remark}
Although the existence of equalizers (kernels) in $\mathbb{S}^{\mathcal{C}}$
is guaranteed, they are \emph{not} necessarily formed in $\mathbb{S}_{A}$
(compare with \cite[Problem 16]{Por2006}). For sufficient conditions for
forming equalizers (kernels) of $\mathcal{C}$-linear maps in $\mathbb{S}%
_{A}, $ see Proposition \ref{Ker-Coker} below.
\end{remark}

The proof of the following result is essentially the same as that for
comodules of corings (\emph{e.g.} \cite{CMZ2002}, \cite{BW2003}).

\begin{proposition}
\label{propert}Let $(\mathcal{C},\Delta _{\mathcal{C}},\varepsilon _{%
\mathcal{C}})$ be an $A$-semicoring and consider the forgetful functor $%
\mathcal{F}:\mathbb{S}^{\mathcal{C}}\longrightarrow \mathbb{S}_{A}.$

\begin{enumerate}
\item For every $M\in \mathbb{S}^{\mathcal{C}},$ we have a functor%
\begin{equation*}
-\otimes _{A}M:\mathbb{S}_{A}\longrightarrow \mathbb{S}^{\mathcal{C}},\text{
}X\mapsto (X\otimes _{A}M,X\otimes _{A}\rho ^{M}).
\end{equation*}%
Moreover, $-\otimes _{A}M$ is left adjoint to $\mathrm{Hom}^{\mathcal{C}%
}(M,-):\mathbb{S}^{\mathcal{C}}\longrightarrow \mathbb{S}_{A};$ we have
natural isomorphisms for all $X_{A}$ and $Y\in \mathbb{S}^{\mathcal{C}}:$%
\begin{equation*}
\mathrm{Hom}^{\mathcal{C}}(X\otimes _{A}M,Y)\simeq \mathrm{Hom}_{A}(X,%
\mathrm{Hom}^{\mathcal{C}}(M,Y)),\text{ }f\longmapsto \lbrack x\longmapsto
\lbrack m\longmapsto f(x\otimes _{A}m)]]
\end{equation*}%
with inverse $g\longmapsto \lbrack x\otimes _{A}m\longmapsto g(x)(m)].$

\item $-\otimes _{A}\mathcal{C}:\mathbb{S}_{A}\longrightarrow \mathbb{S}^{%
\mathcal{C}}$ is right adjoint to $\mathcal{F}.$ We have a natural
isomorphism for all $X_{A}$ and $Y\in \mathbb{S}^{\mathcal{C}}:$%
\begin{equation*}
\mathrm{Hom}^{\mathcal{C}}(Y,X\otimes _{A}\mathcal{C})\simeq \mathrm{Hom}%
_{A}(\mathcal{F}(Y),X),\text{ }f\longmapsto \lbrack y\longmapsto \vartheta
_{X}^{r}\circ (X\otimes _{A}\varepsilon )(f(y))]
\end{equation*}%
with inverse $g\longmapsto \lbrack y\longmapsto \sum g(y_{<0>})\otimes
_{A}y_{<1>}].$
\end{enumerate}
\end{proposition}

\begin{corollary}
\label{-C}Let $(\mathcal{C},\Delta _{\mathcal{C}},\varepsilon _{\mathcal{C}%
}) $ be an $A$-semicoring.

\begin{enumerate}
\item Let $M\in \mathbb{S}^{\mathcal{C}}.$ The functor $-\otimes _{A}M:%
\mathbb{S}_{A}\longrightarrow \mathbb{S}^{\mathcal{C}}$ preserves all
colimits, whence right exact. In particular, it preserves coequalizers \emph{%
(}cokernels\emph{)}, pushouts, \emph{(}regular\emph{)} epimorphisms, direct
sums and directed colimits. On the other hand, the functor $\mathrm{Hom}^{%
\mathcal{C}}(M,-):\mathbb{S}^{\mathcal{C}}\longrightarrow \mathbb{S}_{A}$
preserves all limits, whence left exact. In particular, it preserves
equalizers \emph{(}kernels\emph{)}, pullbacks, monomorphisms, direct
products and inverse limits.

\item $-\otimes _{A}\mathcal{C}:\mathbb{S}_{A}\longrightarrow \mathbb{S}^{%
\mathcal{C}}$ preserves all colimits and all limits, whence exact. In
particular, it preserves coequalizers \emph{(}cokernels\emph{)}, equalizers
\emph{(}kernels\emph{)}, pushouts, pullbacks, \emph{(}regular\emph{)}
epimorphisms, monomorphisms, direct sums, direct products, directed colimits
and inverse limits.

\item The forgetful functor $\mathcal{F}:\mathbb{S}^{\mathcal{C}%
}\longrightarrow \mathbb{S}_{A}$

\begin{enumerate}
\item creates and preserves all colimits, whence right exact. In particular,
it creates and preserves coequalizers \emph{(}cokernels\emph{)}, pushouts,
\emph{(}regular\emph{)} epimorphisms, direct sums and directed colimits.

\item creates all limits which are preserved by $-\otimes _{A}\mathcal{C}:%
\mathbb{S}_{A}\longrightarrow \mathbb{S}_{A}.$
\end{enumerate}
\end{enumerate}
\end{corollary}

The following result provides a sufficient conditions equalizers and kernels
in $\mathbb{S}^{\mathcal{C}}$ to be formed in $\mathbb{S}_{A}.$

\begin{proposition}
\label{Ker-Coker}Let $(\mathcal{C},\Delta ,\varepsilon )$ be $A$-semicorings.

\begin{enumerate}
\item Coequalizers of $\mathbb{S}^{\mathcal{C}}$ are formed in $\mathbb{S}%
_{A}.$ In particular, for any morphism $f:M\longrightarrow N$ in $\mathbb{S}%
^{\mathcal{C}},$ we have $\mathrm{Coker}(f)=N/f(M).$

\item If $_{A}\mathcal{C}$ is flat, then equalizers of $\mathbb{S}^{\mathcal{%
C}}$ are formed in $\mathbb{S}_{A}.$

\item If $_{A}\mathcal{C}$ is $u$-flat, then kernels of $\mathbb{S}^{%
\mathcal{C}}$ are formed in $\mathbb{S}_{A}.$
\end{enumerate}
\end{proposition}

\begin{Beweis}
\begin{enumerate}
\item The forgetful functor $\mathcal{F}$ creates and preserves all
colimits, and in particular coequalizers, by Proposition \ref{Porst}. It
follows that coequalizers (cokernels) are formed in $\mathbb{S}_{A}.$ In
what follows we provide an elementary direct proof. Let $f,g:M%
\longrightarrow N$ be two morphisms in $\mathbb{S}^{\mathcal{C}}$ and let $%
\mathrm{Coeq}(f,g)$ be the coequalizer of $f,g$ in $\mathbb{S}_{A}.$ Since $%
-\otimes _{A}\mathcal{C}:\mathbb{S}_{A}\longrightarrow \mathbb{S}_{A}$
preserves coequalizers, we have the following commutative diagram of right $%
A $-semimodules%
\begin{equation}
\xymatrix{M \ar@<1ex>[rr]^{f} \ar[rr]_{g} \ar[d]_{\rho^M} & & N
\ar[rr]^{\pi} \ar[d]_{\rho^N} & & {\mathrm {Coeq}}(f,g)
\ar@{.>}[d]^{\rho^{{\mathrm {Coeq}}(f,g)}} \\ M\otimes_A {\mathcal C}
\ar@<1ex>[rr]^{g\otimes_A {\mathcal C}} \ar[rr]_{f\otimes_A {\mathcal C}} &
& N \otimes_A {\mathcal C} \ar[rr]^{\pi \otimes_A {\mathcal C}} & & {\mathrm
{Coeq}}(f,g) \otimes_A {\mathcal C} }  \label{coeq}
\end{equation}%
The left square is commutative since $f$ is a morphism of right $\mathcal{C}$%
-semicomodules. It follows that
\begin{eqnarray*}
(\pi \otimes _{A}\mathcal{C})\circ \rho ^{N}\circ f &=&(\pi \otimes _{A}%
\mathcal{C})\circ (f\otimes _{A}\mathcal{C})\circ \rho ^{M} \\
&=&((\pi \circ f)\otimes _{A}\mathcal{C})\circ \rho ^{M} \\
&=&((\pi \circ g)\otimes _{A}\mathcal{C})\circ \rho ^{M} \\
&=&(\pi \otimes _{A}\mathcal{C})\circ (g\otimes _{A}\mathcal{C})\circ \rho
^{M} \\
&=&(\pi \otimes _{A}\mathcal{C})\circ \rho ^{N}\circ g.
\end{eqnarray*}%
By the universal property of coequalizers, there exists a \emph{unique} $A$%
-linear map $\rho ^{{\mathrm{Coeq}}(f,g)}:{\mathrm{Coeq}}(f,g)%
\longrightarrow {\mathrm{Coeq}}(f,g)\otimes _{A}\mathcal{C}$ such that the
right square is commutative. Consider the following diagram with commutative
trapezoids and inner rectangle%
\begin{equation}
\xymatrix{{\mathrm {Coeq}}(f,g) \ar[rrr]^(.45){\rho^{{\mathrm {Coeq}}(f,g)}}
\ar[ddd]_{\rho^{{\mathrm {Coeq}}(f,g)}} & & & {\mathrm {Coeq}}(f,g)
\otimes_A {\mathcal C} \ar[ddd]^{{\mathrm {Coeq}}(f,g) \otimes_A \Delta} \\
{ } & N \ar[d]_(.45){\rho ^N} \ar[r]^(.45){\rho ^N} \ar[ul]^{\pi} & N
\otimes_{A} {\mathcal C} \ar[ur]_(.45){\pi \otimes_A {\mathcal C}}
\ar[d]^(.45){N \otimes_A \Delta} & \\ & N \otimes_{A} {\mathcal {C}}
\ar[dl]^(.45){\pi \otimes_A {\mathcal C}} \ar[r]_(.45){\rho ^N \otimes_A
\mathcal {C}} & N \otimes_{A} {\mathcal {C}} \otimes_{A} {\mathcal {C}}
\ar[dr]^(.45){\pi \otimes_{A} {\mathcal C} \otimes_{A} {\mathcal C}} & \\
{\mathrm {Coeq}}(f,g) \otimes_A {\mathcal C} \ar[rrr]_(.45){\rho^{{\mathrm
{Coeq}}(f,g)} \otimes_A {\mathcal C}} & & & {\mathrm {Coeq}}(f,g) \otimes_A
{\mathcal C} \otimes_A {\mathcal C} }
\end{equation}%
Notice that%
\begin{eqnarray*}
({\mathrm{Coeq}}(f,g)\otimes _{A}\Delta )\circ \rho ^{{\mathrm{Coeq}}%
(f,g)}\circ \pi &=&({\mathrm{Coeq}}(f,g)\otimes _{A}\Delta )\circ (\pi
\otimes _{A}\mathcal{C})\circ \rho ^{N} \\
&=&(\pi \otimes _{A}\mathcal{C}\otimes _{A}\mathcal{C})\circ (N\otimes
_{A}\Delta )\circ \rho ^{N} \\
&=&(\pi \otimes _{A}\mathcal{C}\otimes _{A}\mathcal{C})\circ (\rho
^{N}\otimes _{A}\mathcal{C})\circ \rho ^{N} \\
&=&(\rho ^{{\mathrm{Coeq}}(f,g)}\otimes _{A}\mathcal{C})\circ (\pi \otimes
_{A}\mathcal{C})\circ \rho ^{N} \\
&=&(\rho ^{{\mathrm{Coeq}}(f,g)}\otimes _{A}\mathcal{C})\circ \rho ^{{%
\mathrm{Coeq}}(f,g)}\circ \pi .
\end{eqnarray*}%
Since $\pi $ is an epimorphism, we conclude that%
\begin{equation*}
({\mathrm{Coeq}}(f,g)\otimes _{A}\Delta )\circ \rho ^{{\mathrm{Coeq}}%
(f,g)}=(\rho ^{{\mathrm{Coeq}}(f,g)}\otimes _{A}\mathcal{C})\circ \rho ^{{%
\mathrm{Coeq}}(f,g)}.
\end{equation*}%
Moreover, we have%
\begin{eqnarray*}
\vartheta _{{\mathrm{Coeq}}(f,g)}^{r}\circ ({\mathrm{Coeq}}(f,g)\otimes
_{A}\varepsilon )\circ \rho ^{{\mathrm{Coeq}}(f,g)}\circ \pi &=&\vartheta _{{%
\mathrm{Coeq}}(f,g)}^{r}\circ ({\mathrm{Coeq}}(f,g)\otimes _{A}\varepsilon
)\circ (\pi \otimes _{A}\mathcal{C})\circ \rho ^{N} \\
&=&\pi \circ \vartheta _{N}^{r}\circ (N\otimes _{A}\varepsilon )\circ \rho
^{N} \\
&=&\pi .
\end{eqnarray*}%
Since $\pi $ is an epimorphism, we conclude that $\vartheta _{{\mathrm{Coeq}}%
(f,g)}^{r}\circ ({\mathrm{Coeq}}(f,g)\otimes _{A}\varepsilon )\circ \rho ^{{%
\mathrm{Coeq}}(f,g)}=\mathrm{id}_{{\mathrm{Coeq}}(f,g)}.$ Consequently, $({%
\mathrm{Coeq}}(f,g),\rho ^{{\mathrm{Coeq}}(f,g)})$ is a right $\mathcal{C}$%
-semicomodule.

\item Notice that $\mathbb{S}_{A}$ has equalizers $\mathrm{Eq}(f,g)=\{m\in
M\mid f(m)=g(m)\}$ for any $A$-linear maps $f,g.$ Since $_{A}\mathcal{C}$ is
flat, we have $\mathrm{Eq}(f\otimes _{A}\mathcal{C},g\otimes _{A}\mathcal{C}%
)=\mathrm{Eq}(f,g)\otimes _{A}\mathcal{C}$ in $\mathbb{S}_{A}.$ Consider the
following diagram of right $A$-semimodules%
\begin{equation}
\xymatrix{ {\mathrm {Eq}}(f,g) \ar[rr]^{\iota} \ar@{.>}[d]_{\rho^{{\mathrm
{Eq}}(f,g)}} & & M \ar@<1ex>[rr]^{f} \ar[rr]_{g} \ar[d]_{\rho^M} & & N
\ar[d]^{\rho^ N} \\ {\mathrm {Eq}}(f,g) \otimes_A {\mathcal C}
\ar[rr]^{\iota \otimes_A {\mathcal C}} & & M \otimes_A {\mathcal C}
\ar@<1ex>[rr]^{g \otimes_A {\mathcal C}} \ar[rr]_{f \otimes_A {\mathcal C}}
& & N \otimes_A {\mathcal C} }  \label{eq}
\end{equation}%
Since $f$ is a morphism of right $\mathcal{C}$-semicomodules, the right
square is commutative. It follows that
\begin{eqnarray*}
(f\otimes _{A}\mathcal{C})\circ \rho ^{M}\circ \iota &=&\rho ^{N}\circ
(f\circ \iota ) \\
&=&\rho ^{N}\circ (g\circ \iota ) \\
&=&(g\otimes _{A}\mathcal{C})\circ \rho ^{M}\circ \iota
\end{eqnarray*}%
and so there exists, by the universal property of equalizers, a \emph{unique}
$A$-linear map $\rho ^{\mathrm{Eq}(f,g)}:\mathrm{Eq}(f,g)\longrightarrow
\mathrm{Eq}(f,g)\otimes _{A}\mathcal{C}$ such that the left square is
commutative. Consider the following diagram with commutative trapezoids and
outer rectangle%
\begin{equation}
\xymatrix{M \ar[rrr]^(.45){\rho^{M}} \ar[ddd]_{\rho^{M}} & & & M \otimes_A
{\mathcal C} \ar[ddd]^{M \otimes_A \Delta} \\ { } & {\mathrm {Eq}}(f,g)
\ar[d]_(.45){\rho ^{{\mathrm {Eq}}(f,g)}} \ar[r]^(.45){\rho ^{{\mathrm
{Eq}}(f,g)}} \ar[ul]^{\iota} & {\mathrm {Eq}}(f,g) \otimes_{A} {\mathcal C}
\ar[ur]_(.45){\iota \otimes_A {\mathcal C}} \ar[d]^(.45){{\mathrm {Eq}}(f,g)
\otimes_A \Delta} & \\ & {\mathrm {Eq}}(f,g) \otimes_{A} {\mathcal {C}}
\ar[dl]_(.45){\iota \otimes_A {\mathcal C}} \ar[r]_(.45){\rho ^{{\mathrm
{Eq}}(f,g)} \otimes_A \mathcal {C}} & {\mathrm {Eq}}(f,g) \otimes_{A}
{\mathcal {C}} \otimes_{A} {\mathcal {C}} \ar[dr]^(.45){\iota \otimes_{A}
{\mathcal C} \otimes_{A} {\mathcal C}} & \\ M \otimes_A {\mathcal C}
\ar[rrr]_(.45){\rho^{M} \otimes_A {\mathcal C}} & & & M \otimes_A {\mathcal
C} \otimes_A {\mathcal C} }
\end{equation}%
Notice that%
\begin{eqnarray*}
(\iota \otimes _{A}\mathcal{C}\otimes _{A}\mathcal{C})\circ (\mathrm{Eq}%
(f,g)\otimes _{A}\Delta )\circ \rho ^{\mathrm{Eq}(f,g)} &=&(M\otimes
_{A}\Delta )\circ (\iota \otimes _{A}\mathcal{C})\circ \rho ^{\mathrm{Eq}%
(f,g)} \\
&=&(M\otimes _{A}\Delta )\circ \rho ^{M}\circ \iota \\
&=&(\rho ^{M}\otimes _{A}\mathcal{C})\circ \rho ^{M}\circ \iota \\
&=&(\rho ^{M}\otimes _{A}\mathcal{C})\circ (\iota \otimes _{A}\mathcal{C}%
)\circ \rho ^{\mathrm{Eq}(f,g)} \\
&=&(\iota \otimes _{A}\mathcal{C}\otimes _{A}\mathcal{C})\circ (\rho ^{%
\mathrm{Eq}(f,g)}\otimes _{A}\mathcal{C})\circ \rho ^{\mathrm{Eq}(f,g)}.
\end{eqnarray*}%
Since $_{A}\mathcal{C}$ is flat, it follows that $_{A}\mathcal{C}$ is
mono-flat and so $\iota \otimes _{A}\mathcal{C}\otimes _{A}\mathcal{C}$ is
injective and so%
\begin{equation*}
(\mathrm{Eq}(f,g)\otimes _{A}\Delta )\circ \rho ^{\mathrm{Eq}(f,g)}=(\mathrm{%
Eq}(f,g)\otimes _{A}\mathcal{C})\circ \rho ^{\mathrm{Eq}(f,g)}.
\end{equation*}%
Moreover, we have%
\begin{eqnarray*}
\iota \circ \vartheta _{\mathrm{Eq}(f,g)}^{r}\circ (\mathrm{Eq}(f,g)\otimes
_{A}\varepsilon )\circ \rho ^{\mathrm{Eq}(f,g)} &=&\vartheta _{M}^{r}\circ
(M\otimes _{A}\varepsilon )\circ (\iota \otimes _{A}\mathcal{C})\circ \rho ^{%
\mathrm{Eq}(f,g)} \\
&=&\vartheta _{M}^{r}\circ (M\otimes _{A}\varepsilon )\circ \rho ^{M}\circ
\iota \\
&=&\iota \circ \mathrm{Eq}(f,g).
\end{eqnarray*}%
Since $\mathrm{Eq}(f,g)\overset{\iota }{\hookrightarrow }M$ is a
monomorphism, we conclude that $\vartheta _{\mathrm{Eq}(f,g)}^{r}\circ (%
\mathrm{Eq}(f,g)\otimes _{A}\varepsilon )\circ \rho ^{\mathrm{Eq}(f,g)}=%
\mathrm{id}_{\mathrm{Eq}(f,g)}.$ It follows that $(\mathrm{Eq}(f,g),\rho ^{%
\mathrm{Eq}(f,g)})$ is a right $\mathcal{C}$-semicomodule.

\item Since $_{A}\mathcal{C}$ is $u$-flat, $-\otimes _{A}\mathcal{C}:\mathbb{%
S}_{A}\longrightarrow \mathbb{S}_{A}$ preserves kernels (see Remark \ref%
{kernel-flat}). The proof is along the lines of that of (2).$\blacksquare $
\end{enumerate}
\end{Beweis}

\begin{notation}
Let $\mathcal{C}$ be an $A$-semicoring. In addition to the forgetful functor
$\mathcal{F}:\mathbb{S}^{\mathcal{C}}\longrightarrow \mathbb{S}_{A},$ we
consider the following functors%
\begin{equation*}
\mathcal{G}:=-\otimes _{A}\mathcal{C}:\mathbb{S}_{A}\longrightarrow \mathbb{S%
}^{\mathcal{C}}\text{ and }\mathbb{M}:=-\otimes _{A}\mathcal{C}:\mathbb{S}%
_{A}\longrightarrow \mathbb{S}_{A}.
\end{equation*}
\end{notation}

\begin{remark}
\label{mono-flat}Let $\mathcal{C}$ be an $A$-semicoring. Recall that the
functor $\mathcal{G}$ is exact, whence it preserves monomorphisms and
kernels. If $\mathcal{F}$ preserves monomorphisms (kernels), then $\mathbb{M}%
=\mathcal{F}\circ \mathcal{G}$ preserves monomorphisms (kernels) as well.
\end{remark}

\begin{proposition}
\label{C-U-flat}Let $\mathcal{C}$ be an $A$-semicoring and consider the
forgetful functor $\mathcal{F}:\mathbb{S}^{\mathcal{C}}\longrightarrow
\mathbb{S}_{A}.$

\begin{enumerate}
\item $_{A}\mathcal{C}$ is flat if and only if $\mathcal{F}$ is \emph{(}left%
\emph{)} exact.

\item Assume that $_{A}\mathcal{C}$ is mono-flat. The following are
equivalent:

\begin{enumerate}
\item $_{A}\mathcal{C}$ is uniformly flat;

\item $_{A}\mathcal{C}$ is $u$-flat;

\item $\mathcal{F}$ creates and preserves kernels;

\item $\mathcal{F}$ preserves uniform monomorphisms \emph{(i.e.} every
uniform monomorphism in $\mathbb{S}^{\mathcal{C}}$ is injective\emph{)}.
\end{enumerate}
\end{enumerate}
\end{proposition}

\begin{Beweis}
Recall first that $(\mathcal{F},\mathcal{G})$ and $(\mathcal{G},\mathrm{Hom}%
^{\mathcal{C}}(\mathcal{C},-))$ are adjoint pairs, whence $\mathcal{F}$ is
right exact and $\mathcal{G}$ is exact. Moreover, notice that $\mathbb{M}=%
\mathcal{F}\circ \mathcal{G}:\mathbb{S}_{A}\longrightarrow \mathbb{S}_{A}.$

\begin{enumerate}
\item $(\Rightarrow )$ Notice that we have two (left) exact functors $%
\mathbb{S}_{A}\overset{\mathcal{G}}{\longrightarrow }\mathbb{S}^{\mathcal{C}}%
\overset{\mathcal{F}}{\longrightarrow }\mathbb{S}_{A},$ whence $\mathbb{M}=%
\mathcal{F}\circ \mathcal{G}$ is (left) exact, \emph{i.e. }$_{A}\mathcal{C}$
is flat.

$(\Leftarrow )$ Since $_{A}\mathcal{C}$ is flat, all finite limits are
preserved by the (left) exact functor $\mathcal{G}$ and these are
consequently created and preserved by $\mathcal{F},$ \emph{i.e.} $\mathcal{F}
$ is (left) exact.

\item Assume that $_{A}\mathcal{C}$ is mono-flat.

$(a\Longleftrightarrow b)$ This follows by Remark \ref{kernel-flat}.

$(b\Rightarrow c)$ This follows by Proposition \ref{Ker-Coker}.

$(c\Rightarrow d)$ Let $f:X\longrightarrow Y$ be a uniform monomorphism in $%
\mathbb{S}^{\mathcal{C}}$ and consider the $\ker (f):\mathrm{Ker}%
(f)\longrightarrow X$ in $\mathbb{S}^{\mathcal{C}}.$ Since $f\circ \ker
(f)=0=f\circ 0,$ we conclude that $\mathrm{Ker}(f)=0$ in $\mathbb{S}^{%
\mathcal{C}}.$ Since $\mathcal{F}$ preserves kernels, $\mathrm{Ker}%
(f)=\{0_{X}\}$ in $\mathbb{S}_{A},$ whence $f$ is injective as any $k$%
-uniform $A$-linear map with zero kernel.

$(d\Rightarrow a)$ Let $X\overset{\iota }{\hookrightarrow }Y$ be a uniform $%
A $-subsemimodule, whence $X=\overline{X}=\mathrm{Ker}(Y\overset{\pi }{%
\longrightarrow }Y/X).$ Since $-\otimes _{A}\mathcal{C}:\mathbb{S}%
_{A}\longrightarrow \mathbb{S}^{\mathcal{C}}$ is left exact, it preserves
kernels and so $X\otimes _{A}\mathcal{C}=\mathrm{Ker}(Y\otimes _{A}\mathcal{C%
}\overset{\pi \otimes _{A}\mathcal{C}}{\longrightarrow }Y/X\otimes _{A}%
\mathcal{C})$ in $\mathbb{S}^{\mathcal{C}};$ in particular $\iota \otimes
_{A}\mathcal{C}:X\otimes _{A}\mathcal{C}\longrightarrow Y\otimes _{A}%
\mathcal{C}$ is a uniform monomorphism in $\mathbb{S}^{\mathcal{C}}.$ By our
assumption on $\mathcal{F},$ the map $\iota \otimes _{A}\mathcal{C}:X\otimes
_{A}\mathcal{C}\longrightarrow Y\otimes _{A}\mathcal{C}$ is a uniform
monomorphism in $\mathbb{S}_{A},$ \emph{i.e. }$X\otimes _{A}\mathcal{C}\leq
^{u}Y\otimes _{A}\mathcal{C}.$ We conclude that $_{A}\mathcal{C}$ is
uniformly flat.$\blacksquare $
\end{enumerate}
\end{Beweis}

\begin{remark}
Let $(M,\rho ^{M})$ be a right $\mathcal{C}$-semicomodule and $N\leq _{A}M.$
If $_{A}\mathcal{C}$ is mono-flat and $\rho _{1}^{N},\rho
_{2}^{N}:N\longrightarrow N\otimes _{A}\mathcal{C}$ make $N$ a $\mathcal{C}$%
-subsemicomodule of $M,$ then one can easily see that $\rho _{1}^{N}=\rho
_{2}^{N}.$ However, if $_{A}\mathcal{C}$ is not mono-flat, then it might
happen that $N$ has two \emph{different} structures as a $\mathcal{C}$%
-subsemicomodule of $M.$
\end{remark}

\qquad The following example appeared originally in \cite{Set1974} (and
cited in \cite{Wisch1975}) with $\mathbb{Z}$ at the place of $\mathbb{N}%
_{0}: $

\begin{ex}
Let $C=(\mathbb{N}_{0}\oplus \mathbb{Z}/n\mathbb{Z},\Delta ,\varepsilon )$
be the $\mathbb{N}_{0}$-coalgebra whose comultiplication and counity are
given by%
\begin{eqnarray*}
\Delta ((l,\overline{m})) &:&=(l,0)\otimes _{\mathbb{N}_{0}}(1,0)+(1,0)%
\otimes _{\mathbb{N}_{0}}(0,\overline{m})+(0,\overline{m})\otimes _{\mathbb{N%
}_{0}}(1,0)+(0,\overline{m})\otimes _{\mathbb{N}_{0}}(0,\overline{1}); \\
\varepsilon ((l,\overline{m})) &:&=l.
\end{eqnarray*}%
Consider the Abelian monoid $M=\mathbb{Q}/\mathbb{Z}$ and the embedding of
monoids
\begin{equation*}
\iota :\mathbb{Z}/n\mathbb{Z}\hookrightarrow \mathbb{Q}/\mathbb{Z},\text{ }%
\overline{z}\mapsto \lbrack \frac{r}{n}]\text{ (where }z\equiv r\text{ mod }n%
\text{ and }r\in \{0,1,\cdots ,n-1\}\text{).}
\end{equation*}%
We have a structure of a right $\mathcal{C}$-semicomodule $(M,\rho ^{M})$
and two \emph{different} $\mathcal{C}$-subsemicomodule structures $(N,\rho
_{1}^{N}),$ $(N,\rho _{2}^{N})$ where%
\begin{eqnarray*}
\rho ^{M} &:&\mathbb{Q}/\mathbb{Z}\longrightarrow \mathbb{Q}/\mathbb{Z}%
\otimes _{\mathbb{N}_{0}}C,\text{ }\overline{q}\mapsto \overline{q}\otimes _{%
\mathbb{N}_{0}}(1,0); \\
\rho _{1}^{N} &:&\mathbb{Z}/n\mathbb{Z}\longrightarrow \mathbb{Z}/n\mathbb{Z}%
\otimes _{\mathbb{N}_{0}}C,\text{ }\overline{z}\mapsto \overline{z}\otimes _{%
\mathbb{N}_{0}}(1,0); \\
\rho _{2}^{N} &:&\mathbb{Z}/n\mathbb{Z}\longrightarrow \mathbb{Z}/n\mathbb{Z}%
\otimes _{\mathbb{N}_{0}}C,\text{ }\overline{z}\mapsto \overline{z}\otimes _{%
\mathbb{N}_{0}}(1,0)+(0,\overline{1})\otimes _{\mathbb{N}_{0}}\overline{z}.
\end{eqnarray*}%
Notice that $\mathbb{N}_{0}\oplus \mathbb{Z}/n\mathbb{Z}$ is not mono-flat
in $\mathbf{AbMonoid}.$
\end{ex}

\begin{proposition}
\label{cog}Let $\mathcal{C}$ be an $A$-semicoring.

\begin{enumerate}
\item If $\mathbf{Q}$ is a cogenerator in $\mathbb{S}_{A},$ then $\mathbf{Q}%
\otimes _{A}\mathcal{C}$ is a cogenerator in $\mathbb{S}^{\mathcal{C}}.$

\item If $A_{A}$ is a cogenerator in $\mathbb{S}_{A},$ then $\mathcal{C}$ is
a cogenerator in $\mathbb{S}^{\mathcal{C}}.$
\end{enumerate}
\end{proposition}

\begin{Beweis}
\begin{enumerate}
\item It is well-known that right adjoint functors preserve cogenerators.
However, we provide a direct proof: let $f,g:M\longrightarrow N$ be
morphisms in $\mathbb{S}^{\mathcal{C}}$ with $f\neq g.$ Since $\mathbb{S}%
_{A} $ has products and $\mathbf{Q}_{A}$ is a cogenerator, there exists an
index set $\Lambda $ such that $N\hookrightarrow \mathbf{Q}^{\Lambda }.$ By
Corollary \ref{-C}, $-\otimes _{A}\mathcal{C}$ preserves monomorphisms and
direct products. It follows that we have a monomorphism $\gamma :N\otimes
_{A}\mathcal{C}\longrightarrow \mathbf{Q}^{\Lambda }\otimes _{A}\mathcal{C}%
\simeq (\mathbf{Q}\otimes _{A}\mathcal{C})^{\Lambda }$ in $\mathbb{S}^{%
\mathcal{C}}.$ Consider the canonical projection $\pi :(\mathbf{Q}\otimes
_{A}\mathcal{C})^{\Lambda }\longrightarrow \mathbf{Q}\otimes _{A}\mathcal{C}%
. $ If $\pi \circ \gamma \circ f=\pi \circ \gamma \circ g,$ then we have $%
(\gamma \circ f)(\lambda )=(\gamma \circ g)(\lambda )$ for every $\lambda
\in \Lambda ,$ whence $\gamma \circ f=\gamma \circ g$ and this yields $f=g$
(a contradiction). Setting $h:=\pi \circ \gamma :N\longrightarrow \mathbf{Q}%
\otimes _{A}\mathcal{C},$ we have $h\circ f\neq h\circ g$ and we conclude
that $\mathbf{Q}\otimes _{A}\mathcal{C}$ is a cogenerator in $\mathbb{S}^{%
\mathcal{C}}.$

\item This follows directly from (1) and the canonical isomorphism $A\otimes
_{A}\mathcal{C}\simeq \mathcal{C}.\blacksquare $
\end{enumerate}
\end{Beweis}

\begin{punto}
Let $\mathcal{C}$ be an $A$-semiring with $_{A}\mathcal{C}$ be flat, so that
the forgetful functor $\mathcal{F}:\mathbb{S}^{\mathcal{C}}\longrightarrow
\mathbb{S}_{A}$ is exact by Proposition \ref{C-U-flat} (1). It follows also
that $\mathbb{S}^{\mathcal{C}}$ has kernels (as well as cokernels) formed in
$\mathbb{S}_{A}$ and that monomorphisms are injective while regular
epimorphisms are surjective. One can prove that in this case the category $%
\mathbb{S}^{\mathcal{C}}$ has a $(\mathbf{Surj},\mathbf{Inj})$\emph{%
-factorization system} \cite{AHS2004}. The arguments in \cite{Abu-b} about
the natural definition of exact sequences of semimodules apply to the
category $\mathbb{S}^{\mathcal{C}}$ as well and so we call a sequence$X%
\overset{f}{\longrightarrow }Y\overset{g}{\longrightarrow }Z$ of right $%
\mathcal{C}$-semicomodules \emph{exact }iff $f(X)=\mathrm{Ker}(g)$ and $g$
is $k$-uniform. A sequence $0\longrightarrow X\overset{f}{\longrightarrow }Y%
\overset{g}{\longrightarrow }Z\longrightarrow 0$ will be called a \emph{%
short exact sequence }iff $f$ induces an isomorphism $X\simeq \mathrm{Ker}%
(g) $ and $g$ induces an isomorphism $Z\simeq \mathrm{Coker}(f).$
\end{punto}

\begin{definition}
Let $_{A}\mathcal{C}$ be a semicoring. We say that a right $\mathcal{C}$%
-semicomodule $E$ is \emph{uniformly injective} iff for every uniform
monomorphism $f:M\longrightarrow N$ in $\mathbb{S}^{\mathcal{C}},$ the
induced map of Abelian monoids%
\begin{equation*}
(f,E):\mathrm{Hom}^{\mathcal{C}}(N,E)\longrightarrow \mathrm{Hom}^{\mathcal{C%
}}(M,E),\text{ }h\mapsto h\circ f
\end{equation*}%
is surjective and uniform.
\end{definition}

\begin{remark}
\label{pres-inj}It is well-known that functors between Abelian categories
with an exact left adjoint preserve injective objects \cite[6.28]{Fai1973}.
We extend this result to the functor $\mathcal{G}:=-\otimes _{A}\mathcal{C}:%
\mathbb{S}_{A}\longrightarrow \mathbb{S}^{\mathcal{C}},$ with $_{A}\mathcal{C%
}$ flat, which is right adjoint to the exact forgetful functor. Please
notice that the categories under consideration are, in general, far away
from being Abelian (not even additive).
\end{remark}

\begin{definition}
We say that $\mathbb{S}_{A}$ has \emph{enough }(\emph{uniformly})\emph{\
injective objects} iff every $A$-semimodule is an $A$-subsemimodule of a
(uniformly) injective $A$-semimodule.
\end{definition}

\begin{proposition}
\label{injective}Let $\mathcal{C}$ be an $A$-semicoring and consider the
functor $\mathcal{G}:=-\otimes _{A}\mathcal{C}:\mathbb{S}_{A}\longrightarrow
\mathbb{S}^{\mathcal{C}}.$

\begin{enumerate}
\item If every \emph{(}uniform\emph{)} monomorphism in $\mathbb{S}^{\mathcal{%
C}}$ is injective, then $\mathcal{G}$ preserves \emph{(}uniformly\emph{)}
injective objects.

\item Assume that $\mathbb{S}_{A}$ has enough \emph{(}uniformly\emph{)}
injective objects. Every \emph{(}uniform\emph{)} monomorphism in $\mathbb{S}%
^{\mathcal{C}}$ is injective if and only if $\mathcal{G}$ preserves \emph{(}%
uniformly\emph{)} injective objects.
\end{enumerate}
\end{proposition}

\begin{Beweis}
\begin{enumerate}
\item Let $E$ be a (uniformly) injective $A$-semimodule. Let $\iota
:L\longrightarrow M$ be a (uniform) monomorphism in $\mathbb{S}^{\mathcal{C}%
}.$ By our assumptions $\mathcal{F}$ preserves (uniform) monomorphisms,
whence $L\leq _{A}M$ ($L\leq _{A}^{u}M$). By Proposition \ref{propert}, we
have natural isomorphisms $\mathrm{Hom}^{\mathcal{C}}(L,E\otimes _{A}%
\mathcal{C})\overset{\Phi }{\simeq }\mathrm{Hom}_{A}(\mathcal{F}(L),E)$ and $%
\mathrm{Hom}^{\mathcal{C}}(L,E\otimes _{A}\mathcal{C})\overset{\Phi }{\simeq
}\mathrm{Hom}_{A}(\mathcal{F}(L),E)$ where $\mathcal{F}:\mathbb{S}^{\mathcal{%
C}}\longrightarrow \mathbb{S}_{A}$ is the forgetful functor. Consider the
following commutative diagram of Abelian monoids%
\begin{equation*}
\xymatrix{{\mathrm {Hom}}^{\mathcal{C}}(M,E\otimes _{A}\mathcal{C})
\ar[rr]^{(L,E\otimes _{A}\mathcal{C})} \ar@{=}[d] & & {\mathrm
{Hom}}^{\mathcal{C}}(L,E\otimes _{A}\mathcal{C})\ar@{=}[d] \ar@{.>}[rr] & &
0\\ {\mathrm {Hom}}_A (M,E) \ar[rr]^{(L,E)} & & {\mathrm {Hom}}_A (L,E)
\ar[rr] & & 0}
\end{equation*}

By assumption, $(\iota ,E)$ is a (uniform) surjective map, whence $(\iota
,E\otimes _{A}\mathcal{C})$ is a (uniform) surjective map. Consequently, $%
E\otimes _{A}\mathcal{C}$ is (uniformly) injective in $\mathbb{S}^{\mathcal{C%
}}.$

\item The proof is along the lines of that of the corresponding result for
comodules of coalgebras over a commutative ring \cite[Proposition 8]%
{Wisch1975}. Assume that $\mathcal{G}$ preserves (uniformly) injective
objects. Let $h:L\longrightarrow M$ be a (uniform) monomorphism of right $%
\mathcal{C}$-comodules. We claim that $h$ is injective. By assumption, there
exists a (uniformly) injective right $\mathcal{C}$-semicomodule $E$ such
that $L\overset{\iota _{L}}{\hookrightarrow }E.$ By assumption, $\mathcal{G}$
preserves (uniformly) injective objects, whence $E\otimes _{A}\mathcal{C}$
is (uniformly) injective in $\mathbb{S}^{\mathcal{C}}.$ Notice that we have
a morphisms of $\mathcal{C}$-colinear maps%
\begin{equation*}
(\iota _{L}\otimes _{A}\mathcal{C})\circ \rho ^{L}:L\overset{\rho ^{L}}{%
\longrightarrow }L\otimes _{A}\mathcal{C}\overset{\iota _{L}\otimes _{A}%
\mathcal{C}}{\longrightarrow }E\otimes _{A}\mathcal{C}.
\end{equation*}%
Since $E\otimes _{A}\mathcal{C}$ is (uniformly) injective, there exists a
unique $\mathcal{C}$-colinear map $g:M\longrightarrow E\otimes _{A}\mathcal{C%
}$ such that $g\circ h=(\iota _{L}\otimes _{A}\mathcal{C})\circ \rho ^{L}.$
So, we have%
\begin{eqnarray*}
\iota _{L} &=&\iota _{L}\circ \mathrm{id}_{L} \\
&=&\iota _{L}\circ \vartheta _{L}^{r}\circ (L\otimes _{A}\varepsilon )\circ
\rho ^{L} \\
&=&\vartheta _{E}^{r}\circ (E\otimes _{A}\varepsilon )\circ (\iota
_{L}\otimes _{A}\mathcal{C})\circ \rho ^{L} \\
&=&\vartheta _{E}^{r}\circ (E\otimes _{A}\varepsilon )\circ g\circ h.
\end{eqnarray*}%
It follows that $h$ is injective and we are done.$\blacksquare $
\end{enumerate}
\end{Beweis}

Combining Proposition \ref{C-U-flat} (2) and Proposition \ref{injective}, we
get

\begin{corollary}
\label{ini-flat}Let $\mathcal{C}$ be an $A$-semicoring and consider the
functor $\mathcal{G}:=-\otimes _{A}\mathcal{C}:\mathbb{S}_{A}\longrightarrow
\mathbb{S}^{\mathcal{C}}.$

\begin{enumerate}
\item If $_{A}\mathcal{C}$ is $u$-flat, then $\mathcal{G}$ preserves
uniformly injective objects.

\item Assume that $\mathbb{S}_{A}$ has enough uniformly injective objects
and the $_{A}\mathcal{C}$ is mono-flat. The following are equivalent:

\begin{enumerate}
\item $_{A}\mathcal{C}$ is $u$-flat;

\item every uniform monomorphism in $\mathbb{S}^{\mathcal{C}}$ is injective;

\item $\mathcal{G}$ preserves uniformly injective objects.
\end{enumerate}
\end{enumerate}
\end{corollary}

\begin{proposition}
\label{injec-cog}Let $\mathcal{C}$ be an $A$-semicoring and assume that $_{A}%
\mathcal{C}$ is flat.

\begin{enumerate}
\item If $E_{A}$ is a \emph{(}uniformly\emph{)} injective cogenerator, then $%
E\otimes _{A}\mathcal{C}$ is a \emph{(}uniformly\emph{)} injective
cogenerator in $\mathbb{S}^{\mathcal{C}}.$

\item If $A_{A}$ is a \emph{(}uniformly\emph{)} injective cogenerator, then $%
\mathcal{C}$ is a \emph{(}uniformly\emph{)} injective cogenerator in $%
\mathbb{S}^{\mathcal{C}}.$
\end{enumerate}
\end{proposition}

\section{Measuring $\protect\alpha $-Pairings}

\qquad In this section, we introduce and investigate the $\mathcal{C}$-\emph{%
rational }$\mathcal{A}$\emph{-semimodules} associated with a measuring left
(right) $\alpha $-pairing $(\mathcal{A},\mathcal{C}).$

\subsection*{\textbf{Measuring Pairings}}

\begin{punto}
Let $\mathcal{C}$ be an $A$-semicoring and consider the left dual $A$%
-semiring $^{\ast }\mathcal{C}:=\mathrm{Hom}_{A-}(\mathcal{C},A).$ If $%
\mathcal{A}$ is an $A$-semiring with a morphism of $A$-semirings $\kappa :%
\mathcal{A}\longrightarrow $ $^{\ast }\mathcal{C},$ $a\mapsto \lbrack
c\mapsto <a,c>],$ then we call $P:=(\mathcal{A},\mathcal{C})$ a \emph{%
measuring left }$A$\emph{-pairing}. A \emph{measuring right }$A$\emph{%
-pairing} $P=(\mathcal{A},\mathcal{C})$ consists of an $A$-semiring $%
\mathcal{A}$ and an $A$-semicoring $\mathcal{C}$ with a morphism of $A$%
-semirings $\kappa _{P}:\mathcal{A}\longrightarrow \mathcal{C}^{\ast }.$ If $%
\mathcal{A}$ is an $A$-semiring with a morphism of $A$-semirings $\kappa
_{P}:\mathcal{A}\longrightarrow $ $^{\ast }\mathcal{C}^{\ast },$ then we
call $(\mathcal{A},\mathcal{C})$ a \emph{measuring }$A$\emph{-pairing}.
\end{punto}

\begin{punto}
If $P=(\mathcal{A},\mathcal{C})$ is a measuring left (right) $A$-pairing,
then $\mathcal{C}$ is a right (left) $\mathcal{A}$-semimodule with $\mathcal{%
A}$-action given by%
\begin{equation}
c\leftharpoonup a:=\sum c_{1}<a,c_{2}>\text{ (}a\rightharpoonup c:=\sum
<a,c_{1}>c_{2}\text{).}  \label{A-act}
\end{equation}%
If $P=(\mathcal{A},\mathcal{C})$ is a measuring $A$-pairing, then $\mathcal{C%
}$ is an $(\mathcal{A},\mathcal{A})$-bisemimodule with the right and the
left $\mathcal{A}$-actions in (\ref{A-act}).
\end{punto}

\subsection*{The $\protect\alpha $-Condition}

\begin{punto}
\label{alp-cond}We say that a left $A$-pairing $P=(V,W)\;$\emph{satisfies
the }$\alpha $\emph{-condition}\textbf{, }or is a \emph{left }$\alpha $\emph{%
-pairing} iff the following map is \emph{injective} and \emph{subtractive}
(whence uniform):%
\begin{equation}
\alpha _{M}^{P}:M\otimes _{A}W\ \longrightarrow \mathrm{Hom}_{-A}(V,M),\text{
}\sum m_{i}\otimes _{A}w_{i}\mapsto \lbrack v\mapsto \sum m_{i}<v,w_{i}>].
\notag
\end{equation}%
A right $A$-pairing $P=(V,W)$ is said to \emph{satisfy the }\textbf{(}\emph{%
right}\textbf{) }$\alpha $\emph{-condition}, or to be a\emph{\ right }$%
\alpha $\emph{-pairing,} iff for every left\emph{\ }$A$-semimodule $M,$ the
canonical map $\alpha _{M}^{P}:W\otimes _{A}M\longrightarrow \mathrm{Hom}%
_{A-}(V,M)$ is injective and subtractive.
\end{punto}

\begin{definition}
We say that $_{A}W$ is a \emph{left }$\alpha $\emph{-semimodule} iff the
left $A$-pairing $(^{\ast }W,W)$ satisfies the $\alpha $-condition,
equivalently iff the following canonical map%
\begin{equation*}
\alpha _{M}^{W}:M\otimes _{A}W\longrightarrow \mathrm{Hom}_{-A}(^{\ast }W,M),%
\text{ }m\otimes _{A}w\mapsto \lbrack f\mapsto mf(w)]
\end{equation*}%
is injective and subtractive (uniform). Symmetrically, one defines \emph{%
right }$\alpha $\emph{-semimodules}. Moreover, we say that $_{A}W_{B}$ is an
$\alpha $\emph{-bisemimodule} iff $_{A}W$ and $W_{B}$ are $\alpha $%
-semimodules.
\end{definition}

\begin{remarks}
\begin{enumerate}
\item If $P=(V,W)$ is a left $\alpha $-pairing, then $W\leq _{A}^{u}V^{\ast
} $ (take $M=A$).

\item If $_{A}W$ is finitely projective, then $\mathrm{Ker}(\alpha
_{M}^{W})=0$ for every $M_{A}.$
\end{enumerate}
\end{remarks}

\begin{exs}
A left $A$-semimodule $W$ is an $\alpha $-semimodule if, for example, $W$
satisfies any of the following conditions:

\begin{enumerate}
\item $_{A}W$ is a free $A$-semimodule;

\item $_{A}W$ is a direct summand of a free $A$-semimodule;

\item $_{A}W$ is finitely projective and $\alpha _{M}^{W}$ is subtractive
for every $M_{A}.$
\end{enumerate}
\end{exs}

\begin{lemma}
\label{flat}If $P=(V,W)$ is a measuring left $\alpha $-pairing, then $_{A}W$
is uniformly flat.
\end{lemma}

\begin{Beweis}
Let $M$ be any right $A$-semimodule, $L\leq _{A}^{u}M$ and consider the
commutative diagram of Abelian monoids%
\begin{equation*}
\xymatrix{L \otimes_{A} W \ar[rr]^{\alpha _L ^P} \ar[d]_{\iota _L \otimes_A
W} & & {\rm Hom}_{-{A}}(V,L) \ar@{^{(}->}[d]^{(V,\iota)} \\ M \otimes_{A} W
\ar[rr]_{\alpha_M ^P} & & {\rm Hom}_{-{A}} (V,M)}
\end{equation*}%
It is easy to see that $\mathrm{Hom}_{-A}(V,-)$ preserves uniform morphisms.
By assumption, $\alpha _{L}^{P}$ is injective and uniform, whence $(V,\iota
)\circ \alpha _{L}^{P}$ is injective and, moreover, uniform by Lemma \cite[%
Lemma 1.15]{Abu-b}. It follows that $\alpha _{M}^{P}\circ (\iota _{L}\otimes
_{A}W)=(V,\iota )\circ \alpha _{L}^{P}$ is injective and uniform, whence $%
L\otimes _{A}W\leq _{A}^{u}M\otimes _{A}W.\blacksquare $
\end{Beweis}

The following technical lemma plays an important role in the investigations
of rational semimodules.

\begin{lemma}
\label{q-2}Let $P=(V,W)$ be a left $\alpha $-pairing. If $L$ is a right $A$%
-semimodule and $K\leq _{A}L$ is an $A$-subsemimodule, then we have for
every $\sum l_{i}\otimes _{A}w_{i}\in L\otimes _{A}W:$%
\begin{equation*}
\sum l_{i}\otimes _{A}w_{i}\in \overline{K}\otimes _{A}W\Longleftrightarrow
\sum l_{i}<v,w_{i}>\in \overline{K}\text{ for all }v\in V.
\end{equation*}
\end{lemma}

\begin{Beweis}
Notice that we have an exact sequence of right $A$-semimodules%
\begin{equation*}
0\longrightarrow \overline{K}\overset{\iota _{K}}{\longrightarrow }L\overset{%
\pi _{K}}{\longrightarrow }L/K\longrightarrow 0.
\end{equation*}%
Consider the commutative diagram%
\begin{equation*}
\xymatrix{0 \ar[r] & \overline{K} \otimes_{A} W \ar[rr]^{\iota_K \otimes_A
W} \ar@{^{(}->}[d]^{\alpha_K ^P} & & L \otimes_{A} W \ar[rr]^{\pi_K
\otimes_A W} \ar@{^{(}->}[d]^{\alpha_L ^P} & & L/K \otimes_{A} W
\ar@{^{(}->}[d]^{\alpha _{L/K} ^P} \ar[r] & 0\\ 0 \ar[r] & {\rm Hom}
_{-{A}}(V,\overline{K}) \ar[rr]_{(V,\iota_K)} & & {\rm Hom} _{-{A}}(V,L)
\ar[rr]_{(V,\pi_K)} & & {\rm Hom} _{-{A}}(V,L/K) & }
\end{equation*}%
By Lemma \ref{flat}, $_{A}W$ is uniformly flat and so the first row is
exact. Clearly, $\sum l_{i}<v,w_{i}>\in \overline{K}$ for every $v\in V$ if
and only if $\sum l_{i}\otimes _{A}w_{i}\in \mathrm{Ke}((V,\pi _{K})\circ
\alpha _{L}^{P})=\mathrm{Ke}(\alpha _{L/K}^{P}\circ (\pi _{K}\otimes _{A}W))=%
\mathrm{Ke}(\pi _{K}\otimes _{A}W)=\overline{K}\otimes _{A}W.\blacksquare $
\end{Beweis}

The proof of the following result is similar to that of \cite[Proposition 2.5%
]{AG-TL2001}:

\begin{lemma}
\label{p-2}Let $V,W$ be $(A,A)$-bisemimodules.

\begin{enumerate}
\item If $P=(V,W),$ $P^{\prime }=(V^{\prime },W^{\prime })$ are left $\alpha
$-pairings, then $P\otimes _{A}^{l}P^{\prime }:=(V^{\prime }\otimes
_{A}V,W\otimes _{A}W^{\prime })$ is a left $\alpha $-pairing, where%
\begin{equation*}
\kappa _{P\otimes _{A}^{l}P^{\prime }}(v^{\prime }\otimes _{A}v)(w\otimes
_{A}w^{\prime })=<v,w<v^{\prime },w^{\prime }>>=<<v^{\prime },w^{\prime
}>v,w>.
\end{equation*}

\item If $P=(V,W),$ $P^{\prime }=(V^{\prime },W^{\prime })$ are right $%
\alpha $-pairings, then $P\otimes _{A}^{r}P^{\prime }:=(V\otimes
_{A}V^{\prime },W^{\prime }\otimes _{A}W)$ is a right $\alpha $-pairing,
where
\begin{equation*}
\kappa _{P^{\prime }\otimes _{A}^{r}P}(v\otimes _{A}v^{\prime })(w^{\prime
}\otimes _{A}w)=<v,<v^{\prime },w^{\prime }>w>=<v<v^{\prime },w^{\prime
}>,w>.
\end{equation*}
\end{enumerate}
\end{lemma}

\subsection*{Rational semimodules}

\qquad In what follows, we introduce and investigate the category $\mathrm{%
Rat}^{\mathcal{C}}(\mathbb{S}_{\mathcal{A}})$ of $\mathcal{C}$-\emph{%
rational right }$\mathcal{A}$\emph{-semimodules} associated with a measuring
left $\alpha $-pairing $(\mathcal{A},\mathcal{C}).$

\begin{punto}
\label{rat-dar}Let $P=(\mathcal{A},\mathcal{C})$ a measuring left $\alpha $%
-pairing and $M$ a right $\mathcal{A}$-semimodule. Since $\mathbb{S}_{A}$ is
complete, it has pullbacks. We define $\mathrm{Rat}^{\mathcal{C}}(M_{%
\mathcal{A}})$ as the \emph{pullback} of the following diagram of right $%
\mathcal{A}$-semimodules and $\mathcal{A}$-linear maps
\begin{equation*}
\xymatrix{\mathrm{Rat}^{\mathcal{C}}(M_{\mathcal{A}}) \ar@{^{(}.>}[d]
\ar@{.>}[r]^(.45){\rho ^M} & M \otimes_{A} {\mathcal {C}}
\ar@{^{(}->}[d]^{\alpha_M ^P} \\ M \ar@{^{(}->}[r]_(0.4){\rho_M} & {\rm
Hom}_{-{A}}({\mathcal {A}},M)}
\end{equation*}%
Clearly, $\mathrm{Rat}^{\mathcal{C}}(M_{\mathcal{A}}):=(\rho
_{M})^{-1}(\alpha _{M}^{P}(M\otimes _{A}\mathcal{C})),$ \emph{i.e.} $m\in
\mathrm{Rat}^{\mathcal{C}}(M_{\mathcal{A}})$ iff there exists a uniquely
determined element $\sum m_{i}\otimes _{A}c_{i}\in M\otimes _{A}\mathcal{C}$
such that $ma=\sum m_{i}<a,c_{i}>$ for every $a\in \mathcal{A}.$ We say that
$M_{\mathcal{A}}$ is $\mathcal{C}$\emph{-rational} iff $\mathrm{Rat}^{%
\mathcal{C}}(M_{\mathcal{A}})=M$ and set%
\begin{equation*}
\mathrm{Rat}^{{\mathcal{C}}}(\mathbb{S}_{\mathcal{A}}):=\{M_{\mathcal{A}}%
\text{ }\mid \mathrm{Rat}^{\mathcal{C}}(M_{\mathcal{A}})=M\}.
\end{equation*}%
Symmetrically, if $Q=(\mathcal{A},\mathcal{C})$ is a measuring right $\alpha
$-pairing and $M$ is a left $\mathcal{A}$-semimodule, then we set $^{%
\mathcal{C}}\mathrm{Rat}(_{\mathcal{A}}M):=(\rho _{M})^{-1}(\alpha _{M}^{Q}(%
\mathcal{C}\otimes _{A}M)).$ Similarly, we say that $_{\mathcal{A}}M$ is $%
\mathcal{C}$\emph{-rational} iff $^{\mathcal{C}}\mathrm{Rat}(_{\mathcal{A}%
}M)=M$ and set%
\begin{equation*}
^{{\mathcal{C}}}\mathrm{Rat}(_{\mathcal{A}}\mathbb{S}):=\{\text{ }_{\mathcal{%
A}}M\mid \text{ }^{\mathcal{C}}\mathrm{Rat}(_{\mathcal{A}}M)=M\}.
\end{equation*}
\end{punto}

\begin{punto}
Let $P=(\mathcal{A},\mathcal{C})$ be a measuring left $\alpha $-pairing and $%
Q=(\mathcal{B},\mathcal{D})$ a measuring right $\alpha $-pairing. For each $(%
\mathcal{B},\mathcal{A})$-bisemimodule $(M,\rho _{M}^{\mathcal{A}},\rho
_{M}^{\mathcal{B}}),$ we have%
\begin{equation}
\mathrm{Rat}^{\mathcal{C}}((^{{\normalsize \mathcal{D}}}\mathrm{Rat}(_{%
{\normalsize \mathcal{B}}}M))_{\mathcal{A}})=\text{ }^{{\normalsize \mathcal{%
D}}}\mathrm{Rat}(_{{\normalsize \mathcal{B}}}M){\normalsize \cap }\mathrm{Rat%
}^{\mathcal{C}}(M_{\mathcal{A}})=\text{ }^{{\normalsize \mathcal{D}}}\mathrm{%
Rat}(_{\mathcal{B}}(\mathrm{Rat}^{\mathcal{C}}(M_{\mathcal{A}})))
\label{ratC-D}
\end{equation}%
and set%
\begin{equation*}
^{\mathcal{D}}\mathrm{Rat}^{\mathcal{C}}(_{\mathcal{B}}\mathbb{S}_{\mathcal{A%
}}):=\{_{\mathcal{B}}M_{\mathcal{A}}\mid \mathrm{Rat}^{\mathcal{C}}((^{%
\mathcal{D}}\mathrm{Rat}(_{\mathcal{B}}M))_{\mathcal{A}})=M\}.
\end{equation*}
\end{punto}

The following technical lemma plays an important role in our investigations.

\begin{lemma}
\label{clos}Let $P=(\mathcal{A},\mathcal{C})$ be a measuring left $\alpha $%
-pairing. For every $(M,\rho _{M})\in \mathbb{S}_{\mathcal{A}}$ we have:

\begin{enumerate}
\item $\mathrm{Rat}^{\mathcal{C}}(M_{\mathcal{A}})\subseteq M$ is an $%
\mathcal{A}$-subsemimodule.

\item $\mathrm{Rat}^{\mathcal{C}}(M_{\mathcal{A}})=\overline{\mathrm{Rat}^{%
\mathcal{C}}(M_{\mathcal{A}})}.$

\item For every $L\leq _{\mathcal{A}}M,$ we have $\mathrm{Rat}^{\mathcal{C}}(%
\overline{L}_{\mathcal{A}})=\overline{L}\cap \mathrm{Rat}^{\mathcal{C}}(M_{%
\mathcal{A}}).$

\item $\mathrm{Rat}^{\mathcal{C}}(\mathrm{Rat}^{\mathcal{C}}(M_{\mathcal{A}%
}))=\mathrm{Rat}^{\mathcal{C}}(M_{\mathcal{A}}).$

\item For every $N\in \mathbb{S}_{\mathcal{A}}$ and $f\in $ $\mathrm{Hom}_{-%
\mathcal{A}}(M,N)$ we have $f(\mathrm{Rat}^{\mathcal{C}}(M_{\mathcal{A}%
}))\subseteq \mathrm{Rat}^{\mathcal{C}}(N_{\mathcal{A}}).$
\end{enumerate}
\end{lemma}

\begin{Beweis}
\begin{enumerate}
\item This follows directly from the definition since $\mathbb{S}_{\mathcal{A%
}}$ has pullbacks.

\item This follows from Remark \ref{f-1} (note that $M\otimes _{A}\mathcal{C}%
\hookrightarrow \mathrm{Hom}_{-A}(\mathcal{A},M)$ is subtractive by our
definition of $\alpha $-pairings).

\item Let $m\in \overline{L}\cap \mathrm{Rat}^{\mathcal{C}}(M_{\mathcal{A}})$
with $\rho _{M}(m)=\alpha _{M}^{P}(\sum m_{i}\otimes _{A}c_{i}).$ For every $%
a\in \mathcal{A}$ we have $\sum m_{i}<a,c_{i}>=ma\in \overline{L},$ whence $%
\sum m_{i}\otimes _{A}c_{i}\in \overline{L}\otimes _{A}\mathcal{C}$ by Lemma %
\ref{q-2}, \emph{i.e. }$m\in \mathrm{Rat}^{\mathcal{C}}(\overline{L}).$ The
reverse inclusion is obvious.

\item This follows directly from (2) and (3).

\item Consider the following commutative diagram
\begin{equation*}
\xymatrix{{\mathrm{Rat}}^{\mathcal{C}}(M_{\mathcal{A}})
\ar@{.>}[dr]^(.45){\tilde f} \ar[rr]^(.45){\rho ^M} \ar@{^{(}->}
[dd]_{\iota} & & M \otimes_{A} {\mathcal {C}} \ar[d]^{f \otimes_A {\mathcal
C}} \\ & \mathrm{Rat}^{\mathcal{C}}(N_{\mathcal{A}}) \ar@{^{(}->}[d]
\ar@{->}[r]^(.45){\rho ^N} & N \otimes_{A} {\mathcal {C}}
\ar@{^{(}->}[d]^{\alpha_N ^P} \\ M \ar[r]_{f} & N
\ar@{^{(}->}[r]_(0.4){\rho_N} & {\rm Hom}_{-{A}}({\mathcal {A}},N)}
\end{equation*}%
The equality $\rho _{N}\circ f\circ \iota =(\alpha _{N}^{P}\circ (f\otimes
_{A}\mathcal{C})\circ \rho ^{M}$ and the fact that the inner rectangle is a
pullback, by our definition of $\mathrm{Rat}^{\mathcal{C}}(N_{\mathcal{A}}),$
imply the existence of a unique $A$-linear map $\widetilde{f}:\mathrm{Rat}^{%
\mathcal{C}}(M)\longrightarrow \mathrm{Rat}^{\mathcal{C}}(N)$ which
completes the diagram commutatively. Indeed, $\widetilde{f}=f_{\mid _{%
\mathrm{Rat}^{\mathcal{C}}(M)}}$ and we are done.$\blacksquare $
\end{enumerate}
\end{Beweis}

\begin{remarks}
\label{clos-un}Let $P=(\mathcal{A},\mathcal{C})$ be a measuring left $\alpha
$-pairing.

\begin{enumerate}
\item If $M$ is a $\mathcal{C}$-rational right $\mathcal{A}$-semimodule and $%
L\leq _{\mathcal{A}}^{u}M,$ then it follows from Lemma \ref{clos} (3) that $%
\mathrm{Rat}^{\mathcal{C}}(L_{\mathcal{A}})=L\cap \mathrm{Rat}^{\mathcal{C}%
}(M_{\mathcal{A}})=L\cap M=L,$ \emph{i.e.} $L_{\mathcal{A}}$ is $\mathcal{C}$%
-rational. So, $\mathrm{Rat}^{{\mathcal{C}}}(\mathbb{S}_{\mathcal{A}})$ is
closed under \emph{uniform} subobjects.

\item The embedding $\mathcal{C}\overset{\chi _{P}}{\hookrightarrow }%
\mathcal{A}^{\ast }$ induces an isomorphism of right $\mathcal{A}$%
-semimodules $\mathcal{C}\overset{\chi _{P}}{\simeq }\mathrm{Rat}^{\mathcal{C%
}}((\mathcal{A}^{\ast })_{\mathcal{A}}).$
\end{enumerate}
\end{remarks}

The following results generalize our previous results on rational modules
for corings over associative algebras \cite{Abu2003}.

\begin{proposition}
\label{co-rat}Let $P=(\mathcal{A},\mathcal{C})$ be a left measuring $A$%
-pairing.

\begin{enumerate}
\item We have an embedding%
\begin{equation*}
\iota :\mathbb{S}^{\mathcal{C}}\longrightarrow \mathbb{S}_{\mathcal{A}},%
\text{ }(M,\rho ^{M})\mapsto (M,\alpha _{M}^{P}\circ \rho ^{M}).
\end{equation*}%
In particular, $\mathrm{Hom}^{\mathcal{C}}(M,N)\subseteq \mathrm{Hom}_{%
\mathcal{A}}(M,N)$ for all $M,N\in \mathbb{S}^{\mathcal{C}}.$

\item If $P$ satisfies the $\alpha $-condition, then we have a functor%
\begin{equation*}
\mathrm{Rat}^{\mathcal{C}}(-):\mathbb{S}_{\mathcal{A}}\longrightarrow
\mathbb{S}^{\mathcal{C}},\text{ }(M,\rho _{M})\mapsto (M,(\alpha
_{M}^{P})^{-1}\circ \rho _{M}).
\end{equation*}%
In particular, $\mathrm{Hom}^{\mathcal{C}}(M,N)=\mathrm{Hom}_{\mathcal{A}%
}(M,N)$ for all $M,N\in \mathbb{S}^{\mathcal{C}}.$
\end{enumerate}
\end{proposition}

\begin{proposition}
\label{rat-ad}If $(\mathcal{A},\mathcal{C})$ is a left measuring $\alpha $%
-pairing. The full subcategory $\mathrm{Rat}^{\mathcal{C}}(-)\overset{\iota }%
{\hookrightarrow }\mathbb{S}_{\mathcal{A}}$ is reflective \emph{(i.e.} $%
(\iota ,\mathrm{Rat}^{\mathcal{C}}(-))$ is an adjoint pair of functors\emph{)%
}.
\end{proposition}

We are ready now to present our first main result in this section:

\begin{theorem}
\label{equal}Let $\mathcal{A}$ be an $A$-semiring and $\mathcal{C}$ an $A$%
-semicoring.

\begin{enumerate}
\item If $P=({\mathcal{A}},{\mathcal{C}})$ is a measuring ${}$left $\alpha $%
-pairing, then $\mathbb{S}^{{\mathcal{C}}}\simeq \mathrm{Rat}^{{\mathcal{C}}%
}(\mathbb{S}_{\mathcal{A}}).$

\item If $P=({\mathcal{A}},{\mathcal{C}})$ is a measuring right $\alpha $%
-pairing, then $^{\mathcal{C}}\mathbb{S}\simeq $ $^{\mathcal{C}}\mathrm{Rat}%
(_{\mathcal{A}}\mathbb{S}).$

\item If $P=(\mathcal{A},\mathcal{C})$ is a measuring left $\alpha $-pairing
and $Q=(\mathcal{B},\mathcal{D})$ is a measuring right $\alpha $-pairing,
then $^{\mathcal{D}}\mathbb{S}^{\mathcal{C}}\simeq $ $^{\mathcal{D}}\mathrm{%
Rat}^{\mathcal{C}}(_{\mathcal{B}}\mathbb{S}_{\mathcal{A}}).$
\end{enumerate}
\end{theorem}

\begin{punto}
\label{Biend}For every $A$-semicoring $\mathcal{C}$ we have an isomorphism
of $A$-semirings $(\mathcal{C}^{\ast },\star _{r})\simeq \mathrm{End}^{%
\mathcal{C}}(\mathcal{C})$ via $f\mapsto \lbrack c\mapsto \sum
f(c_{1})c_{2}] $ with inverse $g\mapsto \varepsilon _{\mathcal{C}}\circ g$
(compare with Proposition \ref{propert} (1)). Symmetrically, $(^{\ast }%
\mathcal{C},\star _{l})\simeq $ $^{\mathcal{C}}\mathrm{End}(\mathcal{C}%
)^{op} $ as $A$-semirings. If $P=(\mathcal{A},\mathcal{C})$ is a measuring
left $\alpha $-pairing, then we have by Proposition \ref{co-rat} (2) $%
\mathcal{C}^{\ast }\simeq \mathrm{End}^{\mathcal{C}}(\mathcal{C})=\mathrm{End%
}(\mathcal{C}_{\mathcal{A}}).$ On the other hand, if $P$ is a measuring
right $\alpha $-pairing, then $^{\ast }\mathcal{C}\simeq $ $^{\mathcal{C}}%
\mathrm{End}(\mathcal{C})^{op}=\mathrm{End}(_{\mathcal{A}}\mathcal{C})^{op}.$
In particular, if $\mathcal{C}$ satisfies the left and the right $\alpha $%
-conditions then we have%
\begin{equation*}
\mathrm{End}(_{\mathcal{C}^{\ast }}\mathcal{C}_{^{\ast }\mathcal{C}})=\text{
}^{\mathcal{C}}\mathrm{End}^{\mathcal{C}}(\mathcal{C})\simeq Z(\mathcal{C}%
^{\ast })=Z(^{\ast }\mathcal{C}).
\end{equation*}
\end{punto}

An important role by studying the category of rational representations
related to a left measuring $\alpha $-pairings is played by the following
finiteness results which holds for the restricted class of completely
subtractive semicomodules.

\begin{lemma}
\label{M-fin}Let $P=(\mathcal{A},\mathcal{C})$ be a measuring left $\alpha $%
-pairing. If $M\in \mathrm{Rat}^{\mathcal{C}}(\mathbb{S}_{A})$ is completely
subtractive, then there exists for every finite subset $\{m_{1},...,m_{k}\}%
\subset M$ some $N\in \mathrm{Rat}^{\mathcal{C}}(\mathbb{S}_{A})$, such that
$N\subset M$ and $N_{A}$ is finitely generated.
\end{lemma}

\begin{Beweis}
Let $\{m_{1},\cdots ,m_{k}\}\subset M.$ For each $i=1,\cdots ,n,$ we have $%
m_{i}\mathcal{A}\leq _{\mathcal{A}}^{u}M,$ whence a $\mathcal{C}$%
-subsemicomodule by Remark \ref{clos-un} and Proposition \ref{co-rat}.
Moreover $m_{i}\in m_{i}\mathcal{A}=\overline{m_{i}\mathcal{A}}$ and
consequently there exists a subset $\{(m_{ij},c_{ij})\}_{j=1}^{n_{i}}\subset
m_{i}\mathcal{A}\times \mathcal{C}$ such that $\rho
_{M}(m_{i})=\sum\limits_{j=1}^{n_{i}}m_{ij}\otimes _{A}c_{ij}$ for $%
i=1,...,k.$ Obviously, $N:=\sum\limits_{i=1}^{k}m_{i}\mathcal{A}%
=\sum\limits_{i=1}^{k}\sum\limits_{j=1}^{n_{i}}m_{ij}A\leq M$ is a $\mathcal{%
C}$-subsemicomodule and contains $\{m_{1},...,m_{k}\}.\blacksquare $
\end{Beweis}

An application of Lemma \ref{M-fin} and its dual yields the following
finiteness result:

\begin{proposition}
\label{C-fin}Let $\mathcal{C}$ be an $A$-semicoring. If $\mathcal{C}_{A}$
\emph{(}resp. $_{A}\mathcal{C},$ $_{A}\mathcal{C}_{A}$\emph{)} is completely
subtractive, then every finite subset of $\mathcal{C}$ is contained in a
right $\mathcal{C}$-coideal \emph{(}resp. left $\mathcal{C}$-coideal, $%
\mathcal{C}$-bicoideal\emph{)}, which is finitely generated in $\mathbb{S}%
_{A}.$
\end{proposition}

\begin{lemma}
\label{subgen}Let $\mathcal{C}$ be an $A$-semicoring.

\begin{enumerate}
\item Every right $\mathcal{C}$-semicomodule is a subsemicomodule of a $%
\mathcal{C}$-generated right $\mathcal{C}$-semicomodule.

\item $\mathbb{S}^{\mathcal{C}}\subseteq \sigma \lbrack \mathcal{C}_{^{\ast }%
\mathcal{C}}].$
\end{enumerate}
\end{lemma}

\begin{Beweis}
\begin{enumerate}
\item Let $(M,\rho ^{M})$ be an arbitrary right $\mathcal{C}$-semicomodule.
There exists a set $\Lambda $ and a surjective morphism of right $A$%
-semimodules $A^{(\Lambda )}\overset{\pi }{\longrightarrow }M\longrightarrow
0$ \cite[Proposition 17.11]{Gol1999}. It follows that we have a surjective
morphism right $\mathcal{C}$-semicomodules $C^{(\Lambda )}\simeq A^{(\Lambda
)}\otimes _{A}\mathcal{C}\overset{\pi \otimes _{A}\mathcal{C}}{%
\longrightarrow }M\otimes _{A}\mathcal{C}\longrightarrow 0.$ So, $M\otimes
_{A}\mathcal{C}$ is generated by $\mathcal{C}$ as an object of $\mathbb{S}^{%
\mathcal{C}}.$ Since $\rho ^{M}:M\longrightarrow M\otimes _{A}\mathcal{C}$
is a retraction and $\mathcal{C}$-colinear, we conclude that $M$ is a
subobject of $M\otimes _{A}\mathcal{C}$ in $\mathbb{S}^{\mathcal{C}}.$

\item Since morphisms of right $\mathcal{C}$-semicomodules are $^{\ast }%
\mathcal{C}$-linear by Proposition \ref{co-rat}, the results follows by (1).$%
\blacksquare $
\end{enumerate}
\end{Beweis}

\begin{lemma}
\label{alp-inj}Let $(\mathcal{C},\Delta ,\varepsilon )$ be an $A$-semicoring
such that $\alpha _{M}^{\mathcal{C}}:M\otimes _{A}\mathcal{C}\longrightarrow
\mathrm{Hom}_{-A}(^{\ast }\mathcal{C},M)$ is subtractive for every $M_{A}.$
If $\mathbb{S}^{\mathcal{C}}=\sigma _{u}[\mathcal{C}_{^{\ast }\mathcal{C}}],$
then $_{A}\mathcal{C}$ is a mono-flat $\alpha $-semimodule.
\end{lemma}

\begin{Beweis}
Assume that $\mathbb{S}^{\mathcal{C}}=\sigma _{u}[\mathcal{C}_{^{\ast }%
\mathcal{C}}].$ Notice that in this case every monomorphism in $\mathbb{S}^{%
\mathcal{C}}$ is injective, whence $_{A}\mathcal{C}$ is mono-flat by Remark %
\ref{mono-flat}. Let $M$ be an arbitrary right $A$-semimodule and consider $%
(M\otimes _{A}\mathcal{C},M\otimes _{A}\Delta )\in \mathbb{S}^{\mathcal{C}%
}=\sigma _{u}[\mathcal{C}_{^{\ast }\mathcal{C}}].$ For every $L\in \sigma
_{u}[\mathcal{C}_{^{\ast }\mathcal{C}}],$ we have a commutative diagram%
\begin{equation*}
\xymatrix{ {\mathrm {Hom}}^{\mathcal C}(L,M\otimes_A {\mathcal C})
\ar@{=}[rr] \ar[d]_{\simeq} & & {\mathrm {Hom}}_{-{^{\ast }\mathcal{C}}}
(L,M\otimes_A {\mathcal C}) \ar[d]^{(L,\alpha_M ^{\mathcal C})}\\ {\mathrm
{Hom}}_{-A} (L,M) \ar[d]_{\simeq} & & {\mathrm {Hom}}_{-{^{\ast
}\mathcal{C}}} (L,{\mathrm {Im}}(\alpha_M ^{\mathcal C})) \ar@{^{(}->}[d] \\
{\mathrm {Hom}}_{-A}(L \otimes_{{^{\ast }\mathcal{C}}} {^{\ast
}\mathcal{C}}, M) \ar[rr]_{\simeq} & & {\mathrm {Hom}}_{-{^{\ast
}\mathcal{C}}} (L,{\mathrm {Hom}}_{-A} ({^{\ast }\mathcal{C}},M))}
\end{equation*}%
and so $(L,\alpha _{M}^{\mathcal{C}})$ is injective. It follows that $\alpha
_{M}^{\mathcal{C}}:M\otimes _{A}\mathcal{C}\longrightarrow \func{Im}(\alpha
_{M}^{\mathcal{C}})$ is a monomorphism in $\sigma _{u}[\mathcal{C}_{^{\ast }%
\mathcal{C}}],$ whence injective.$\blacksquare $
\end{Beweis}

We are now ready to present the second main result in this section.

\begin{theorem}
\label{eq-alpha}Let $(\mathcal{C},\Delta ,\varepsilon )$ be an $A$%
-semicoring such that $\alpha _{M}^{\mathcal{C}}$ is subtractive for every $%
M_{A}.$

\begin{enumerate}
\item The following are equivalent:

\begin{enumerate}
\item $\mathbb{S}^{\mathcal{C}}=\sigma \lbrack \mathcal{C}_{^{\ast }\mathcal{%
C}}];$

\item $_{A}\mathcal{C}$ is a mono-flat $\alpha $-semimodule and $\mathbb{S}^{%
\mathcal{C}}$ is closed under $^{\ast }\mathcal{C}$-subsemimodules.

\item $\mathbb{S}^{\mathcal{C}}$ is a full subcategory of $\mathbb{S}%
_{^{\ast }\mathcal{C}}$ and is closed under $^{\ast }\mathcal{C}$%
-subsemimodules.

In this case we have%
\begin{equation}
\mathbb{S}^{\mathcal{C}}\simeq \mathrm{Rat}^{\mathcal{C}}(\mathbb{S}_{^{\ast
}\mathcal{C}})=\sigma _{u}[\mathcal{C}_{^{\ast }\mathcal{C}}]=\sigma \lbrack
\mathcal{C}_{^{\ast }\mathcal{C}}].  \label{sub=u-sub}
\end{equation}
\end{enumerate}

\item If $\mathcal{C}$ is a left subtractive $A$-semicoring, then the
following are equivalent:

\begin{enumerate}
\item $\mathbb{S}^{\mathcal{C}}=\sigma _{u}[\mathcal{C}_{^{\ast }\mathcal{C}%
}];$

\item $_{A}\mathcal{C}$ is a mono-flat $\alpha $-semimodule.

\item $\mathbb{S}^{\mathcal{C}}$ is a full subcategory of $\mathbb{S}%
_{^{\ast }\mathcal{C}}.$
\end{enumerate}
\end{enumerate}
\end{theorem}

\begin{Beweis}
\begin{enumerate}
\item $(1)\Rightarrow (2)$ By Lemma \ref{alp-inj}, $_{A}\mathcal{C}$ is a
mono-flat $\alpha $-semimodule. Moreover, $\sigma \lbrack \mathcal{C}%
_{^{\ast }\mathcal{C}}]$ is -- by definition -- closed under $^{\ast }%
\mathcal{C}$-subsemimodules.

$(2)\Rightarrow (3)$ By Theorem \ref{equal} (1), $\mathbb{S}^{\mathcal{C}%
}\simeq \mathrm{Rat}^{\mathcal{C}}(\mathbb{S}_{^{\ast }\mathcal{C}%
})\hookrightarrow \mathbb{S}_{^{\ast }\mathcal{C}}$ is a full subcategory.

$(3)\Rightarrow (1)$ Since $\mathbb{S}^{\mathcal{C}}$ is cocomplete and
closed under homomorphic images, it follows that $\mathrm{Gen}(\mathcal{C}%
_{^{\ast }\mathcal{C}})\subseteq \mathbb{S}^{\mathcal{C}}\subseteq \sigma
\lbrack \mathcal{C}_{^{\ast }\mathcal{C}}]$ (the last inclusion follows by
Lemma \ref{subgen}). However, $\sigma \lbrack \mathcal{C}_{^{\ast }\mathcal{C%
}}]$ is -- by definition -- the smallest subclass of $\mathbb{S}_{^{\ast }%
\mathcal{C}}$ which contains $\mathrm{Gen}(\mathcal{C}_{^{\ast }\mathcal{C}%
}) $ and is closed under $^{\ast }\mathcal{C}$-subsemimodules, whence $%
\mathbb{S}^{\mathcal{C}}=\sigma \lbrack \mathcal{C}_{^{\ast }\mathcal{C}}].$

\item We need only to prove $(3)\Rightarrow (1):$ As in (1), we have $%
\mathrm{Gen}(\mathcal{C}_{^{\ast }\mathcal{C}})\subseteq \mathbb{S}^{%
\mathcal{C}}\subseteq \sigma _{u}[\mathcal{C}_{^{\ast }\mathcal{C}}],$ where
the last inclusions follows by Lemma \ref{subgen} and our assumptions on the
$A$-semicoring which imply that $M\leq _{A}^{u}M\otimes _{A}\mathcal{C}$ for
each $M\in \mathbb{S}^{\mathcal{C}}.$ Notice also that $\mathbb{S}^{\mathcal{%
C}}\simeq \mathrm{Rat}^{\mathcal{C}}(\mathbb{S}_{^{\ast }\mathcal{C}})$ is
closed under \emph{uniform} $^{\ast }\mathcal{C}$-subsemimodules by Remark %
\ref{clos-un} (1), whence $\mathbb{S}^{\mathcal{C}}=\sigma _{u}[\mathcal{C}%
_{^{\ast }\mathcal{C}}]$ since $\sigma _{u}[\mathcal{C}_{^{\ast }\mathcal{C}%
}]$ is -- by definition -- the smallest subcategory of $\mathbb{S}_{^{\ast }%
\mathcal{C}}$ which contains $\mathrm{Gen}(\mathcal{C}_{^{\ast }\mathcal{C}%
}) $ and is closed under uniform $^{\ast }\mathcal{C}$-subsemimodules.$%
\blacksquare $
\end{enumerate}
\end{Beweis}

\begin{proposition}
\label{MC=AM-coring}Let $\mathcal{C}$ be an $A$-semicoring and consider the
functors $\mathbf{R}:=-\otimes _{A}\mathcal{C}:\mathbb{S}_{A}\longrightarrow
\mathbb{S}_{^{\ast }\mathcal{C}}$ and $\mathbb{M}:=-\otimes _{A}\mathcal{C}:%
\mathbb{S}_{A}\longrightarrow \mathbb{S}_{A}.$ If $\mathbb{S}^{\mathcal{C}}=%
\mathbb{S}_{^{\ast }\mathcal{C}},$ then

\begin{enumerate}
\item $\#_{\eta _{^{\ast }\mathcal{C}}}\simeq \mathcal{F}$ and $\mathbf{R}%
\simeq -\otimes _{A}$ $^{\ast }\mathcal{C};$

\item $\mathbf{R}$ has a \emph{(left)} exact adjoint $\mathbf{L}$ such that $%
\mathbb{M}=\mathbf{L}\circ \mathbf{R};$

\item $_{A}\mathcal{C}$ is flat;

\item The forgetful functor $\mathcal{F}:\mathbb{S}^{\mathcal{C}%
}\longrightarrow \mathbb{S}_{A}$ is \emph{(}left\emph{)} exact;

\item $\mathbb{M}:=-\otimes _{A}\mathcal{C}:\mathbb{S}_{A}\longrightarrow
\mathbb{S}_{A}$ is (left) exact;

\item If $_{A}\mathcal{C}$ is uniformly generated and $A^{\Lambda }\otimes
_{A}-:$ $_{A}\mathbb{S}\longrightarrow \mathbf{AbMonoid}$ preserves $i$%
-uniform morphisms, then:

\begin{enumerate}
\item $_{A}\mathcal{C}$ is uniformly finitely presented.

\item $_{A}\mathcal{C}$ is finitely presented.

\item $_{A}\mathcal{C}$ is finitely generated and projective.
\end{enumerate}
\end{enumerate}
\end{proposition}

\begin{Beweis}
\begin{enumerate}
\item The morphism of $A$-semirings $\eta _{^{\ast }\mathcal{C}%
}:A\longrightarrow $ $^{\ast }\mathcal{C}$ induces an adjoint pair of
functors $(-\otimes _{A}$ $^{\ast }\mathcal{C},\#_{\eta _{^{\ast }\mathcal{C}%
}}),$ where $\#_{\eta _{^{\ast }\mathcal{C}}}:\mathbb{S}_{^{\ast }\mathcal{C}%
}\longrightarrow \mathbb{S}_{A}$ is the so called \emph{restriction of
scalars functor}. Indeed, $\#_{\eta _{^{\ast }\mathcal{C}}}\simeq \mathcal{F}
$ in our case and so we have $\mathbf{R}\simeq -\otimes _{A}$ $^{\ast }%
\mathcal{C}$ by the uniqueness of the left adjoint functor.

\item Notice that $\mathcal{F}$ is (left) exact and is left adjoint to $%
\mathbf{R}\simeq \mathcal{G}:=-\otimes _{A}\mathcal{C}:\mathbb{S}%
_{A}\longrightarrow \mathbb{S}^{\mathcal{C}}.$

\item Since $\mathbf{R}$ has a left adjoint, it preserves equalizers, whence
$\mathbb{M}=\mathbf{L}\circ \mathbf{R}$ preserves equalizers and it follows
that $_{A}\mathcal{C}$ is flat.

\item this follows by Proposition \ref{C-U-flat} (1).

\item Notice that both $\mathcal{G}:=-\otimes _{A}\mathcal{C}:\mathbb{S}%
_{A}\longrightarrow \mathbb{S}^{\mathcal{C}}$ and the forgetful functor $%
\mathcal{F}:S^{\mathcal{C}}\longrightarrow S_{A}$ are (left) exact, whence $%
\mathbb{M}=\mathcal{F}\circ \mathcal{G}:\mathbb{S}_{A}\longrightarrow
\mathbb{S}_{A}$ is (left) exact.

\item Since $\mathbb{M}$ is left exact, it preserves products. In
particular, $A^{\Lambda }\otimes _{A}\mathcal{C}\simeq \mathcal{C}^{\Lambda
} $ for every index set $\Lambda $ and it follows that $_{A}\mathcal{C}$ is
uniformly finitely presented by Lemma \ref{f.p.} (2) (notice that we assumed
that $A^{\Lambda }\otimes _{A}-$ preserves $i$-uniform morphisms). Since $%
\mathbb{M}$ is left exact, $_{A}\mathcal{C}$ is flat (by definition).
Finally, as indicated in Lemma \ref{f.p.} (5), finitely presented flat
semimodules are projective.$\blacksquare $
\end{enumerate}
\end{Beweis}

\begin{definition}
We say that an $A$-semicoring $\mathcal{C}$ is \emph{left uniform} iff $_{A}%
\mathcal{C}$ is uniformly generated, $\alpha _{M}^{\mathcal{C}}:M\otimes _{A}%
\mathcal{C}\longrightarrow \mathrm{Hom}_{-A}(^{\ast }\mathcal{C},M)$ is
subtractive for every $M_{A}$ and every $\mathcal{C}$-generated right $%
^{\ast }\mathcal{C}$-semimodules is completely subtractive. Symmetrically,
we define \emph{right uniform }$A$-semicorings. A left and right subtractive
$A$-semicoring is said to be \emph{uniform}.
\end{definition}

\begin{theorem}
\label{c=s}Let $\mathcal{C}$ be a left subtractive $A$-semicoring and assume
that $A^{\Lambda }\otimes _{A}-$ preserves $i$-uniform morphisms. The
following are equivalent:

\begin{enumerate}
\item $\mathbb{S}^{\mathcal{C}}=\mathbb{S}_{^{\ast }\mathcal{C}};$

\item $\#_{\eta _{^{\ast }\mathcal{C}}}\simeq \mathcal{F}$ and $\mathbf{R}%
\simeq -\otimes _{A}$ $^{\ast }\mathcal{C};$

\item $\mathbf{R}$ has a \emph{(left)} exact adjoint $\mathbf{L}$ such that $%
\mathbb{M}=\mathbf{L}\circ \mathbf{R}:\mathbb{S}_{A}\longrightarrow \mathbb{S%
}_{A};$

\item $_{A}\mathcal{C}$ is flat and \emph{(}uniformly\emph{)} finitely
presented;

\item $_{A}\mathcal{C}$ is finitely generated and \emph{(}finitely\emph{)}
projective;

\item $_{\mathcal{C}^{\ast }}\mathcal{C}$ is finitely generated and $_{A}%
\mathcal{C}$ is an $\alpha $-semimodule.
\end{enumerate}
\end{theorem}

\begin{Beweis}
The implications $(1)\Rightarrow (2)\Rightarrow (3)\Rightarrow
(4)\Rightarrow (5)$ are clear by following the proof of Proposition \ref%
{MC=AM-coring}, which was given in this order.

$(5)\Rightarrow (6)$ Since $_{A}\mathcal{C}$ is uniformly finitely generated
and $\mathcal{C}^{\ast }\otimes _{A}-$ preserves cokernels, whence normal
quotients, it follows that $_{^{\ast }\mathcal{C}}\mathcal{C}$ is uniformly
finitely generated. Moreover, since $_{A}\mathcal{C}$ is finitely
projective, we have $\mathrm{Ker}(\alpha _{M}^{\mathcal{C}})=0$ for each $%
M_{A},$ whence $\alpha _{M}^{\mathcal{C}}$ is injective (notice that $\alpha
_{M}^{\mathcal{C}}$ is assumed to be subtractive).

$(6)\Rightarrow (1)$ Notice that $\mathcal{C}$ is a faithful and finitely
generated as a left $\mathcal{C}^{\ast }$-semimodule. Since $\mathcal{C}%
^{\ast }\simeq \mathrm{End}^{\mathcal{C}}(\mathcal{C})=\mathrm{End}(\mathcal{%
C}_{^{\ast }\mathcal{C}}),$ it follows by and that
\begin{equation*}
\mathbb{S}^{\mathcal{C}}\overset{\text{Theorem \ref{eq-alpha} (2)}}{=}\sigma
_{u}[\mathcal{C}_{^{\ast }\mathcal{C}}]\overset{\mathcal{C}\text{ is left
subtractive}}{=}\sigma \lbrack \mathcal{C}_{^{\ast }\mathcal{C}}]\overset{%
\text{Proposition \ref{sg=all}}}{=}\mathbb{S}_{^{\ast }\mathcal{C}%
}.\blacksquare
\end{equation*}
\end{Beweis}

\begin{theorem}
\label{summary}Let $\mathcal{C}$ be an $A$-semicoring such that $_{A}%
\mathcal{C}$ is uniformly generated, $\alpha _{M}^{\mathcal{C}}$ is
subtractive for every $M_{A}$ and assume that $A^{\Lambda }\otimes _{A}-$
preserves $i$-uniform morphisms. We have $\mathbb{S}^{\mathcal{C}}=\mathbb{S}%
_{^{\ast }\mathcal{C}}$ if and only if $_{A}\mathcal{C}$ is finitely
generated and projective and $\mathbb{S}^{\mathcal{C}}$ is closed under $%
^{\ast }\mathcal{C}$-subsemimodules.
\end{theorem}

\begin{Beweis}
$(\Rightarrow )$ Follows by Proposition \ref{MC=AM-coring} and Lemma \ref%
{f.p.}

$(\Longleftarrow )$ As in the proof of Theorem \ref{c=s}, we have $_{%
\mathcal{C}^{\ast }}\mathcal{C}$ is finitely generated and $_{A}\mathcal{C}$
is an $\alpha $-semimodule, whence%
\begin{equation*}
\mathbb{S}^{\mathcal{C}}\overset{\text{Theorem \ref{eq-alpha} (1)}}{=}\sigma
\lbrack _{^{\ast }\mathcal{C}}\mathcal{C}]\overset{\text{Proposition \ref%
{sg=all}}}{=}\mathbb{S}_{^{\ast }\mathcal{C}}.\blacksquare
\end{equation*}
\end{Beweis}

\end{document}